\documentclass[12pt, leqno]{article}

\usepackage{amsmath,amssymb,amsfonts, amscd, theorem, latexsym, graphicx}
\usepackage[all]{xy}

\newtheorem{theo}{Theorem}[section]
\newtheorem{lem}[theo]{Lemma}
\newtheorem{cor}[theo]{Corollary}
\newtheorem{prop}[theo]{Proposition}
\newtheorem{defi}[theo]{Definition}

{\theorembodyfont{\rmfamily} \newtheorem{rem}[theo]{Remark}}
{\theorembodyfont{\rmfamily} \newtheorem{ex}[theo]{Example}}
{\theorembodyfont{\rmfamily} }
{\theorembodyfont{\rmfamily} \newtheorem{con}[theo]{Conjecture}}

\DeclareFontFamily{U}{rsf}{}
\DeclareFontShape{U}{rsf}{m}{n}{
  <5> <6> rsfs5 <7> <8> <9> rsfs7 <10->  rsfs10}{}
\DeclareMathAlphabet{\mathscr}{U}{rsf}{m}{n}
\newcommand{\mycal}[1]{\mathscr{#1}}

\newlength{\unten}
\newlength{\gesamt}
\DeclareFontFamily{U}{cyr}{}
\DeclareFontShape{U}{cyr}{m}{n}{
  <5> wncyr5 <6> wncyr6 <7> wncyr7 <8> wncyr8 <9> wncyr9 <10->
wncyr10}{}
\DeclareMathAlphabet{\mathcyr}{U}{cyr}{m}{n}

\input cyracc.def
\input epsf

\newcommand{\op}[1]{\operatorname{#1}}
\newcommand{\TSh}{{\mathcyr{\cyracc Sh}}}

\newcommand{\sff}[1]{{\sf{#1}}}

\newcommand{\gp}[6]{
\xymatrix@R+10pt@C+15pt{
#1 \ar[d] \ar@<.5ex>[drr]^-{s} \ar@<-.5ex>[drr]_-{t} & & & \\
{#2}\times_{#3} {#2} \ar@<.5ex>[rr]^-{#4} \ar@<-.5ex>[rr]_-{#5} & & #2
\ar[r]^-{#6} & #3
}
}

\newcommand{\bF}{\boldsymbol{F}}
\newcommand{\bS}{\boldsymbol{S}}

\newcommand{\btheta}{\boldsymbol{\theta}}
\newcommand{\bT}{\boldsymbol{T}}
\newcommand{\boldeta}{\boldsymbol{\eta}}
\newcommand{\bzeta}{\boldsymbol{\zeta}}
\newcommand{\bupsilon}{\boldsymbol{\upsilon}}

\newtheorem{theorem}{Theorem}[section]

\theoremstyle{definition}

\theoremstyle{remark}
\newtheorem{remark}[theorem]{Remark}

\numberwithin{equation}{section}

\begin{document}
\title{Twisting Derived Equivalences} 
\author{Oren Ben-Bassat}

\maketitle

\begin{abstract}
We introduce a new method for ``twisting'' relative equivalences of
derived categories of sheaves on two spaces over the same base.  The
first aspect of this is that the derived categories of
sheaves on the spaces are twisted.  They become derived categories of
sheaves on gerbes living over spaces that are locally (on the base)
isomorphic to the original spaces.  Secondly, this is done in a
compatible way so that the equivalence is maintained.  We apply this
method by proving the conjectures of Donagi and Pantev on dualities
between gerbes on genus-one fibrations and comment on other
applications to families of higher genus curves.  We also include a
related conjecture in Mirror Symmetry.  This represents a modified
version of my May 2006 thesis for the University of Pennsylvania.  The
only addition is the subsection titled 'Alternative Method'.
\end{abstract}

\tableofcontents


\section{Introduction}
Categories of sheaves of modules on geometric objects seem to play an
important role in algebraic geometry.  It is sometimes possible to
analyze the algebraic properties of such categories in detail, which
in turn sheds light on the nature of the space itself.  On the other
hand, they also suggest a broader perspective in which a space can
vary in ``non-geometric'' directions, and eventually the notion of a
space could be replaced by a category with certain properties.
Gabriel showed in \cite{Gabriel} that one can recover a Noetherian
scheme from its category of coherent sheaves.  Alternatively, one can
study the derived category of coherent sheaves.  This is a less rigid
structure, and certainly allows for different spaces to have
equivalent derived categories.  This ``derived equivalence'' game
started with Mukai's equivalence of dual complex tori \cite{mukai} and
we give a quick summary in what follows.  The fundamental question
which motivates this work (although we only scratch its surface), is,
{\it How does the derived category vary in families?}.  Although the
derived category is certainly known not to glue, in some cases it
seems to behave as if it does.  That is to say, if one is careful
enough, one can prove things that would easily follow if descent for
derived categories held.  In this thesis, we investigate the special
case when the gluing takes place on the level of the abelian category
of sheaves (which as a stack can be glued from its restriction to a
cover).  Rather than directly study these questions for the derived
category of a space, we study the variance (or twisting) of pairs of
derived categories related by a Fourier-Mukai type equivalence.  In
practice, we start with a given Fourier-Mukai equivalence between the
derived categories of two spaces $X$ and $Y$.  We ask: If the category
of sheaves on $X$ is twisted in some way, can the derived category of
$Y$ be twisted in a compatible way in order to recover a new
equivalence?  On the infinitesimal level, this question was addressed
by Y. Toda in \cite{Toda}.  A formal analysis to all orders was
completed, in a special case, in \cite{BBP}.  A complete picture of
dualities for arbitrary gerby deformation-quantizations of abelian
schemes has been completed by D. Arinkin, building on the results from
his 2002 Thesis \cite{DimaThesis}.  A non-formal version of this same
case was carried out by J. Block in \cite{Block1, Block2}.
Topological versions of twisted T-dualities have been studied by
Mathai and Rosenberg in \cite{MR1,MR2} and by Bunke, Rumpf, and Schick
in \cite{BRS}.  There is also groupoid approach due to these dualities
due to C. Daenzer \cite{Calder}.  It is somewhat interesting and
strange that the philosophy of Toda from \cite{Toda} will be helpful
to us, even though our twistings are not necessarily deformations, and
certainly not formal deformations.  According to this philosophy, for
every derived equivalence between $X$ and $Y$, one should be able to
first find an algebraic object describing the twistings of the derived
category $\sff{D}(X)$ compatible with the derived equivalence.  Next
one should find a natural isomorphism with these twistings of
$\sff{D}(X)$ and the analogous twistings of $\sff{D}(Y)$.  Finally,
one should find an equivalence between the categories associated to a
pair of compatible twistings.  In this thesis, we synthesize two types
of twistings.  The first are those coming from replacing a space by
another space locally isomorphic to it.  The second are those coming
from replacing the derived category of sheaves on a space with the
derived category of sheaves on a gerbe over the space.  The most
convenient setting therefore becomes the one where $X$ and $Y$ are
both fibered over the same space $B$, and the derived equivalence
respects that structure.  In this case, our twists of the derived
categories correspond to decomposing and then re-gluing the relative
stack of the abelian categories of sheaves of the two spaces.  We can
now state our main theorem (with slight re-wording), to be proven as
\ref{theo:main}.  We use here the definitions of $\Phi$-compatible and
$\Phi$-dual which can be found as definitions \ref{defi:compatable}
and \ref{defi:dual} respectively.

\begin{theorem}
Let $X$ and $Y$ be compact, connected, complex manifolds, mapping to a
complex analytic space $B$, via maps $\pi: X \to B$ and $\rho:Y \to
B$, where $\rho$ is flat.  Let $\mathcal{P}$ be a coherent sheaf, flat
over $Y$, on the fiber product $X \times_{B} Y$, which gives an
equivalence of categories $\Phi: \sff{D}^{b}_{c}(Y) \to
\sff{D}^{b}_{c}(X)$, $\Phi = R \phi$ with
\[\phi(\mathcal{S}) = \widetilde{\rho}_{*}(\mathcal{P} \otimes 
\widetilde{\pi}^{*}\mathcal{S}).
\]
Then for any $\Phi$-compatible gerbe $\mathfrak{X}$ over a twisted
version of $X \to B$, and any $\Phi$-dual gerbe $\mathfrak{Y}$ to
$\mathfrak{X}$, there is an equivalence of categories 
\[\widetilde{\Phi}: \sff{D}^{b}_{c}(\mathfrak{Y}, -1) \to 
\sff{D}^{b}_{c}(\mathfrak{X}, -1),
\]
where $\widetilde{\Phi} = R \widetilde{\phi}$ and $\widetilde{\phi}$
is locally built out of $\phi$.
\end{theorem}

The motivation for such a theorem begins at least with Dolgachev and
Gross \cite{DG}, who relate the Tate-Shafarevich group of an elliptic
fibration to its Brauer group.  This naturally leads one to ask if the
derived category of a geometrically twisted elliptic fibration (a
genus one fibration) is somehow related to a derived category of
sheaves on a gerbe over the elliptic fibration.  Such a relationship
was indeed shown in \cite{Andrei1, Andrei2}.  More general dualities
involving gerbes on genus one fibrations were proven in \cite{DP}.
The conjectures made in \cite{DP} form the main geometric motivation
for this thesis, although the methods of proof, will be closer to
those found in \cite{Andrei2}.  More recently, these kind of results
have been used and expanded upon by other authors, for example see
\cite{BriMor}, \cite{BK}, and \cite{Sawon}.  

After the thesis was submitted, we added an alternative, and in many
ways stronger, version of the main theorem which can be found as
\ref{theo:main2}.

In the last section, we comment on some future, more exotic
applications of the ideas we have developed including applications to
homological mirror symmetry, and hyperelliptic families in algebraic
geometry.

\section{Notation and Conventions}
\subsection{Basic Information}
In this thesis, we work with complex manifolds, or complex analytic
spaces and we always use the classical (analytic) topology.  If
$\mathcal{S}$ is a sheaf of abelian groups on a topological space $X$,
and $\mathfrak{U}$ is a cover of $X$, 
then we denote by $\check{C}_{\mathfrak{U}}(X,\mathcal{S})$,
$\check{Z}_{\mathfrak{U}}(X,\mathcal{S})$, and
$\check{B}_{\mathfrak{U}}(X,\mathcal{S})$ the \v{C}ech co-chains,
co-cycles, and co-boundaries, for the sheaf of abelian groups
$\mathcal{S}$ computed with respect to the cover $\mathfrak{U}$ of
$X$.  When we are dealing with a ringed space (or stack), by a {\it
sheaf} (with no qualifications) we always mean a sheaf of modules for
the structure sheaf of the space (or stack).

For our purposes, it will always suffice to consider stacks on a
complex analytic space $X$ which are defined on the site whose
underlying category consists of analytic spaces over $X$, and whose
covering families for a given $W \to X$ consist of surjective local
isomorphisms $\{W' \to W\}$ over $X$.  We will often consider the case
$W=X$.  For a map $W \to X$ we often use the convention 
\[W^{k} = W \times_{X} W \times_{X} W \times_{X} \cdots \times_{X} 
W   \ \ \ \ \ \ \ \ \ (\text{the}\ \ (k+1)-\text{fold product})
\]
Since our methods are quite general, we expect that the
results contained in this thesis apply in other contexts such as the
\'{e}tale topology on schemes, or for twisting different kinds of
relative dualities in mathematics.  

By a gerbe, we shall always mean a $\mathcal{O}^{\times}$-gerbe with
trivial band on a complex analytic space.  On a complex analytic space
$X$, we use $\mathfrak{Mod}(X)$ to denote the stack of abelian
categories of $\mathcal{O}_{X}$-modules, which has global sections
$\text{Mod}(X)$.  We use $\mathfrak{Coh}(X)$ ($\mathfrak{QCoh}(X)$) to
denote the stack of abelian categories of coherent (quasi-coherent)
$\mathcal{O}_{X}$-modules, which has global sections the abelian
category $\text{Coh}(X)$ ($\text{QCoh}(X)$).  We will also need to
consider sheaves on gerbes $\mathfrak{X} \to X$.  For every analytic
space $Z$ mapping to $X$ by $f: Z \to X$, we consider the stack of
functors $f^{-1}\mathfrak{X} \to f^{-1} \mathfrak{Mod}(X)$ and their
natural transformations.  This gives the stack of abelian categories
of $\mathcal{O}_{\mathfrak{X}}$-modules $\mathfrak{Mod}(\mathfrak{X})
\to X$, and similarly we have the stack of abelian categories of
coherent (quasi-coherent) $\mathcal{O}_{\mathfrak{X}}$-modules
$\mathfrak{Coh}(\mathfrak{X}) \to X$ ($\mathfrak{QCoh}(\mathfrak{X})
\to X$) and we call their global sections the abelian categories
$\text{Mod}(\mathfrak{X})$ and $\text{Coh}(\mathfrak{X})$
($\text{QCoh}(\mathfrak{X})$) respectively.  We use $\sff{D}^{*}(X)$
and $\sff{D}^{*}(\mathfrak{X})$ to denote the derived categories of
$\text{Mod}(X)$ and $\text{Mod}(\mathfrak{X})$.  When nothing appears
in the location of the symbol $*$, or $*=\emptyset$, this refers to the
unbounded derived categories.  When $*=-$, this refers to the bounded
above derived category, and when $*=b$, this refers to the bounded
derived category.  Also $\sff{D}_{c}^{*}(X)$ and
$\sff{D}_{c}^{*}(\mathfrak{X})$ or $\sff{D}_{qc}^{*}(X)$ and
$\sff{D}_{qc}^{*}(\mathfrak{X})$ refer to the derived categories of
coherent and quasi-coherent sheaves.  For a discussion of
quasi-coherence in the complex analytic context, see \cite{BBP}.  The
sheaves of weight $k$ on a gerbe will be denoted by
$\mathfrak{Mod}(\mathfrak{X},k)$, $\mathfrak{Coh}(\mathfrak{X},k)$,
$\mathfrak{QCoh}(\mathfrak{X},k)$ with global sections
$\text{Mod}(\mathfrak{X},k)$, $\text{Coh}(\mathfrak{X},k)$,
$\text{QCoh}(\mathfrak{X},k)$, with their associated derived
categories $\sff{D}^{*}(\mathfrak{X},k)$,
$\sff{D}_{c}^{*}(\mathfrak{X},k)$, and
$\sff{D}_{qc}^{*}(\mathfrak{X},k)$.


\section{Basic Techniques}
\subsection{Mukai's Insight} \label{ssec:MukaiYoga}  

In this section we review some well known facts, mainly due to Mukai
\cite{mukai}.  None of the material in this section is original.

Given two complex analytic spaces $M$ and $N$ over $B$ and an element
\[
K\in\text{Mod}(M\times_{B}N),
\]
the integral transform
\[
\phi_{K}^{[M\to N]} : \text{Mod}(M)\to 
\text{Mod}(N)
\]
is defined by
\[
\phi_{K}^{[M\to N]}(G) = p_{N*}( p_{M}^{*}G \otimes K).
\]
Also for a map 
$f:M \to N$ over $B$, we have 
\[\phi_{\mathcal{O}_{\Gamma_{f}}}^{[M\to N]} \cong f_{\star}
\]
as functors $\text{Mod}(M) \to \text{Mod}(N)$
and 
\[\phi_{\mathcal{O}_{\Gamma_{f}}}^{[N\to M]} \cong f^{\star}
\]
as functors $\text{Mod}(N) \to
\text{Mod}(M)$.  For $L$ a line bundle on a space $M$,
we use $T_{L}:\text{Mod}(M) \to
\text{Mod}(M)$ to denote the functor 
\[T_{L}(\mathcal{M}) := L \otimes \mathcal{M} \hspace{20 mm}
\forall \mathcal{M} \in \text{ob}(\text{Mod}(M))
\]
\[T_{L}(h) := \text{id}_{L} \otimes h  \hspace{20 mm}
\forall h \in \text{Hom}_{M}(\mathcal{M}', \mathcal{M}'')
\]

Then we have 
\[ \phi_{\Delta_{*}L}^{[M\to M]} \cong T_{L}
\]
where $\Delta: M \to M \times_{B} M$ is the diagonal embedding.  This
isomorphism uses the symmetry property of the tensor product.  Similarly,
for $f$ an automorphism of $M$, and $L$ a line bundle on $M$, we can
implement the functor $T_{L} \circ f^{*}$ by the object $(f,1)_{*}L
\in \text{Mod}(M \times_{B} M)$.  In other words, we
have an isomorphism
\begin{equation} \label{eqn:represent}
T_{L} \circ f^{*} \cong \phi_{(f,1)_{*}L}^{[M\to M]}
\end{equation}
This isomorphism uses the symmetry property of the tensor product. 

Given two complex analytic spaces $M$ and $N$ and an element $\mathcal{K} \in
\sff{D}(M \times_{B} N)$ we denote by $\Phi_{K}^{[M \to N]}$ the functor
\[\sff{D}(M) \to \sff{D}(N)
\]
given by
\[ \mathcal{S} \mapsto R\pi_{N*}(\mathcal{K} 
\otimes^{L} L\pi_{M}^{*}\mathcal{S} ).
\]

Now, as noted in \cite{AH}, these transforms enjoy the following
locality property.  Fix flat maps $h:S \to B$ and $h':S' \to S$.  For
us the most important case will be the inclusions of an open sets.  Let
$M_{S}, N_{S}, M_{S'}, N_{S'}$ denote the fiber products.  Let
$\mathcal{K}_{S}$, and $K_{S'}$ denote the derived pullbacks of
$\mathcal{K}$ to $\sff{D}(M_{S} \times_{S} N_{S})$ and $\sff{D}(M_{S'}
\times_{S'} N_{S'})$ respectively.  Let $h'_{M}: M_{S'} \to M_{S}$ and
$h'_{N}: N_{S'} \to N_{S}$ be the canonical morphisms Then there is a
natural isomorphism of functors
\begin{equation} \label{eqn:locality}
L {h'}_{N}^{*} \circ \Phi_{\mathcal{K}_{S}} 
\cong \Phi_{\mathcal{K}_{S'}} \circ Lh_{M}^{*}
\end{equation}

The integral transform has the following convolution property (see
\cite{mukai} or \cite[Proposition~11.1]{polishchuk}): If $M$, $N$ and
$P$ are complex analytic spaces over $B$ and $\mathcal{K} \in
\sff{D}^{b}(M \times_{B} N)$ and $\mathcal{L} \in \sff{D}^{b}(N
\times_{B} P)$, then one has a natural isomorphism of functors
\[
\Phi_{\mathcal{L}}^{[N \to P]} \circ 
\Phi_{\mathcal{K}}^{[M \to N]} 
\cong \Phi_{\mathcal{L}*\mathcal{K}}^{[M \to P]},
\]
where
\[
\mathcal{L}*\mathcal{K} =R{p_{M \times_{B} 
P}}_{*}(Lp^{*}_{N \times_{B} 
P}\mathcal{L} \otimes^{L}
  L p^{*}_{M \times_{B} N}\mathcal{K} ) 
\quad \in \sff{D}^{b}(M \times_{B} P),
\]
and $p_{M \times N}$, $p_{M\times 
P}$, $p_{N \times P}$ are the natural
projections $M \times_{B} N \times_{B} 
P \to M \times_{B} N$, etc.

This isomorphism of functors is compatible with base change, in the
sense that if $\mathcal{L}_{S}$ and $\mathcal{K}_{S}$ are the derived
pullbacks of $\mathcal{L}$ and $\mathcal{K}$ inside $\sff{D}(M_{S}
\times_{S} N_{S})$ and \\ $\sff{D}(N_{S} \times_{S} P_{S})$, with respect
to the maps $h_{MN}: M_{S} \times_{S} N_{S} \to M \times_{B} N$ and \\
$h_{NP}: N_{S} \times_{S} P_{S} \to N \times_{B} P$ and $h_{MP}: M_{S}
\times_{S} P_{S} \to M \times_{B} P$ induced from $h: S \to B$.
 
\[(Lh_{NP}^{*}\mathcal{L}) \star (Lh_{MN}^{*}\mathcal{K}) 
\cong Lh_{MP}^{*}(\mathcal{L} \star \mathcal{K}) \in 
\sff{D}(M_{S} \times_{S} P_{S})
\]


Also, from Huybrechts' book \cite{Huybrechts}, we have that
\begin{itemize}
\item
for any $f:N \to P$ over $B$ we have the natural isomorphism of functors
$\sff{D}^{b}(M) \to \sff{D}^{b}(P)$
\begin{equation} \label{eqn:Huy1}
f_{*} \circ \Phi_{K}^{[M\to N]} \cong 
\Phi_{(1_{M} \times f)_{*}K}^{[M\to P]} 
\end{equation}
\item for $f: P \to N$ over $B$ we have the natural isomorphism of functors
$\sff{D}^{b}(M) \to \sff{D}^{b}(P)$
\begin{equation} \label{eqn:Huy2}
f^{*} \circ \Phi_{K}^{[M\to N]} \cong 
\Phi_{(1_{M} \times f)^{*}K}^{[M\to P]} 
\end{equation}
\item for $g:W \to M$ over $B$ we have the natural isomorphism of functors
$\sff{D}^{b}(M) \to \sff{D}^{b}(P)$
\begin{equation} \label{eqn:Huy3}
\Phi_{K}^{[M\to N]} \circ g_{*} \cong 
\Phi_{(g \times 1)^{*}K}^{[W\to N]}
\end{equation}
\item for $g: M \to W$ over $B$ we have the natural isomorphism of
functors $\sff{D}^{b}(M) \to \sff{D}^{b}(P)$
\begin{equation} \label{eqn:Huy4}
\Phi_{K}^{[M\to N]} \circ g^{*} \cong 
\Phi_{(g \times 1)_{*}K}^{[W\to N]}
\end{equation}
\end{itemize}
Of course we
also have the analogue of equation \ref{eqn:represent} in the derived
context, which is an isomorphism of functors $\sff{D}^{b}(M) \to
\sff{D}^{b}(M)$, for any automorphism $f$ of $M$, and line bundle $L$
on $M$
\begin{equation} \label{eqn:derivedrepresent}
T_{L} \circ f^{*} \cong \Phi_{(f,1)_{*}L}^{[M\to M]}
\end{equation}


For a gerbe $\mathfrak{X}$, it is a generally known fact that the
category $\text{Mod}(\mathfrak{X})$ has enough injective and enough
flat objects, see for example \cite{Andrei2} and \cite{lieblich}.
Indeed, following \cite{lieblich}, if $\mathfrak{a}: 
\mathcal{U} \to \mathfrak{X}$
is an atlas, and $\mathcal{M}$ is a sheaf on $\mathfrak{X}$, then we
can find an injection
\[
\mathfrak{a}^{*}\mathcal{M} \to \mathcal{I}
\]
where $\mathcal{I}$ is an injective 
$\mathcal{O}_{\mathcal{U}}$-module.  We can 
also find 
a surjection
\[\mathcal{F} \to \mathfrak{a}^{*}\mathcal{M}.
\]
where $\mathcal{F}$ is a flat $\mathcal{O}_{\mathcal{U}}$-module.
Therefore, we have a sequence of injections
\begin{equation} \label{eqn:eqn-inj_inj} \mathcal{M} 
\to \mathfrak{a}_{*}\mathfrak{a}^{*}\mathcal{M} 
\to \mathfrak{a}_{*}\mathcal{I},
\end{equation}
where $\mathfrak{a}_{*}\mathcal{I}$ is injective and a sequence of surjections
\begin{equation} \label{eqn:eqn-flat_proj} \mathfrak{a}_{!}\mathcal{F} 
\to \mathfrak{a}_{!}\mathfrak{a}^{*}\mathcal{M} 
\to \mathcal{M}
\end{equation}
where $\mathfrak{a}_{!}\mathcal{F}$ is flat.

\subsection{Some Functors}
Recall that for $X$ and $Y$ complex analytic spaces, a gerbe
$\mathfrak{Y} \to Y$ and a morphism $f: X \to Y$ we can form the
$f^{-1} \mathcal{O}_{Y}^{\times}$ banded gerbe $f^{-1} \mathfrak{Y}$,
and then via the map of sheaves of groups $f^{-1}
\mathcal{O}_{Y}^{\times} \to \mathcal{O}_{X}^{\times}$, we define an
($\mathcal{O}_{X}^{\times}$-banded) gerbe $f^{*} \mathfrak{Y}$ on $X$.

We will often use the left
exact functor $f_{*}$, and the right exact functor $f^{*}$
\begin{itemize}
\item $f_{*}: \text{Mod}(f^{*}\mathfrak{Y}) \to
\text{Mod}(\mathfrak{Y})$
\item $f^{*}: \text{Mod}(\mathfrak{Y}) \to
\text{Mod}(f^{*}\mathfrak{Y})$
\end{itemize}
The first of these is defined by the composition
\[\mathfrak{Y} \to f_{*}f^{*} \mathfrak{Y} \to f_{*} 
\mathfrak{Mod}(X) \to 
\mathfrak{Mod}(Y).
\]
The second is defined by the composition 
\[f^{*} \mathfrak{Y} \to f^{*} \mathfrak{Mod}(Y) 
\to \mathfrak{Mod}(X).
\] 

Also if we have two gerbes $\mathfrak{X} \to X$ and 
$\mathfrak{Y} \to X$ over the same space we can define the gerbe 
$\mathfrak{X} \otimes \mathfrak{Y} \to X$ to be the 
$\mathcal{O}_{X}^{\times}$ gerbe 
induced from the $\mathcal{O}_{X}^{\times} 
\times \mathcal{O}_{X}^{\times}$ gerbe 
$\mathfrak{X} \times \mathfrak{Y}$ via the map of sheaves of groups 
$\mathcal{O}_{X}^{\times} \times 
\mathcal{O}_{X}^{\times} \to \mathcal{O}_{X}^{\times}$, $(a,b) \mapsto ab$. 

If we are given $\mathcal{M} \in
\op{ob}(\text{Mod}(\mathfrak{X}))$ and 
$\mathcal{N} \in \op{ob}(\text{Mod}(\mathfrak{Y}))$, then we can form
\[
\mathcal{M} \otimes \mathcal{N} \in
\op{ob}(\text{Mod}(\mathfrak{X} \otimes \mathfrak{Y} ))
\]
This leads, for $\mathcal{N} \in 
\op{ob}(\text{Mod}(\mathfrak{X}))$  to the right exact 
functor $T_{\mathcal{N}}$
\[
T_{\mathcal{N}} : \text{Mod}(\mathfrak{X}) 
\to \text{Mod}(\mathfrak{X}\otimes\mathfrak{Y}).
\]
This is given by the composition 
\[\mathfrak{X} \otimes \mathfrak{Y} \to \mathfrak{X} \otimes 
\mathfrak{Mod}(X) \to \mathfrak{Mod}(X)
\otimes \mathfrak{Mod}(X) \to
\mathfrak{Mod}(X)
\]
\subsection{Some Derived Functors}

A morphism $f:X \to Y$ is {\em flat} if the stalk $\mathcal{O}_{X,x}$
is a flat $\mathcal{O}_{Y, f(x)}$ module for all points $x \in X$.
The flatness of $f$ is equivalent to the exactness of $f^{*}$.  If
$\mathcal{F}$ is a sheaf on $X$, then $\mathcal{F}$ is said to be {\em
flat over} $Y$ if $\mathcal{F}_{x}$ is a flat $\mathcal{O}_{Y, f(x)}$
module for all points $x \in X$.  If $\mathcal{F}$ is flat over $Y$
and $f$ is flat, then the functor $T_{\mathcal{F}} \circ f^{*}$ is exact.

If $f$ is flat, we have the derived functor
\[
f^{*}:
\sff{D}(\mathfrak{Y}) \to
\sff{D}(f^{*}\mathfrak{Y}),
\] obtained by applying
$f^{*}$ term by term to complexes.  This functor obviously preserves
the three types of boundedness.  For a general map $f$ we have the
derived functors
\[
Rf_{*}:\sff{D}^{+}(f^{*}\mathfrak{Y}) \to
\sff{D}^{+}(\mathfrak{Y})
\]
\[
Lf^{*}:\sff{D}^{-}(\mathfrak{Y}) \to
\sff{D}^{-}(f^{*}\mathfrak{Y})
\]
and 
\[
LT_{\mathcal{N}} : 
\sff{D}^{-}(\mathfrak{X}) 
\to \sff{D}^{-}(\mathfrak{X}\otimes \mathfrak{Y}). 
\]

\subsection{Presentations of Gerbes} \label{ssec:GerbePresent}
\medskip
Let $W$ be an analytic space.  A presentation of a gerbe on $W$ is defined 
to be the following.
\begin{itemize}
\item
an analytic space $\mathcal{U}$ mapping to $W$ by a surjective local
isomorphism (atlas) to $W$:
\[a: \mathcal{U} \to W,
\]
\item
a line bundle $L$ over $\mathcal{U} \times_{W}
\mathcal{U}$; we will denote the left and right 
$\mathcal{O}-$module structures on $L$ by
\[l:\mathcal{O} \otimes_{\mathbb{C}} L \to L \ \ \ \ \text{and} 
\ \ \ \ r: L \otimes_{\mathbb{C}} \mathcal{O} \to L \] 
\item 
an isomorphism
\[\theta: p_{12}^{*}L \otimes p_{01}^{*}L \to p_{02}^{*}L \]
over $\mathcal{U} \times_{W} \mathcal{U} \times_{W} \mathcal{U}$
\item
an isomorphism 
\[
\eta:\mathcal{O} \to \Delta^{*}L \]
over $\mathcal{U}$, where $\Delta:\mathcal{U} \to 
\mathcal{U} \times_{W} \mathcal{U}$ is the diagonal map
\end{itemize}
satisfying the conditions

\begin{itemize}
\item that there is an equality of isomorphisms 
\[p_{013}^{*}\theta \circ (p_{123}^{*} \theta \times 1_{p_{01}^{*}L}) = 
p_{023}^{*} \theta \circ (1_{p_{23}^{*}L} \times p_{012}^{*} \theta)\]
\item
\[(1_{\mathcal{U}} \times \Delta)^{*} 
\theta \circ (p_{1}^{*} \eta \times 1_{L}) = r\]
\item
\[(\Delta \times 1_{\mathcal{U}})^{*} 
\theta \circ (1_{L} \times p_{0}^{*} \eta) = l\]
\end{itemize}
We will denote the presentation of a gerbe over $W$ defined in this way by
\[{}_{(L, \theta, \eta)} W.
\]
Here, we suppress the information of the map $a: \mathcal{U} \to W$ in
the notation, assuming that this map is clear from the context.

\begin{defi} \label{defi:strong_equiv_of_pres}
A strong equivalence of two presentations of gerbes with atlas $a:
\mathcal{U} \to W$ ${}_{(L, \theta, \eta)} W$ and ${}_{(L', \theta',
\eta')} W$ is an isomorphism $\mu: L \to L'$ that satisfies
\[
p_{02}^{*} \mu \circ
\theta = \theta' \circ ((p_{12}^{*} \mu) \otimes (p_{01}^{*} \mu))
\]
and
\[\delta^{*}\mu \circ \eta = \eta'.
\]
\end{defi}

For completeness we note that a weak equivalence of two presentations
of gerbes ${}_{(L, \theta, \eta)} W$ and ${}_{(L', \theta', \eta')} W$
over $a: \mathcal{U} \to W$ is a pair $(Q, \tau)$ consisting of
\begin{itemize}
\item
A line bundle $Q \to \mathcal{U}$.
\item
An isomorphism
\[\tau: L \otimes p_{0}^{*}Q \to p_{1}^{*}Q \otimes L'
\]
\end{itemize}
Satisfying the equality of isomorphisms
\[p_{12}^{*}L \otimes p_{01}^{*}L \otimes p_{0}^{*}Q \to 
p_{2}^{*}Q \otimes p_{02}^{*}L'\]
given by 
\[p_{02}^{*} \tau  \circ (\theta \otimes 1_{p_{0}^{*}Q}) = 
(1_{p_{2}^{*}Q} \otimes \theta') \circ 
(p_{12}^{*} \tau \otimes 1_{p_{01}^{*}L'}) 
\circ (1_{p_{12}^{*}L} \otimes p_{01}^{*}\tau)
\]
and the equality of isomorphisms
\[
Q\otimes \Delta^{*}L \to Q
\]
given by 
\[r \circ (1_{Q} \otimes \eta)^{-1} = l \circ 
(\eta'\otimes 1_{Q})^{-1} \circ \Delta^{*}\tau
\]
where we have inserted the natural isomorphisms 
\[Q \otimes \Delta^{*}L \cong \Delta^{*}(p_{1}^{*}Q \otimes L) \ \ \ \
\text{and} \ \ \ \ \Delta^{*}L' \otimes M 
\cong \Delta^{*}(L' \otimes p_{0}^{*}Q).\]


\subsection{Sheaves on Presentations of Gerbes} \label{ssec:TwistedSheaves}

A sheaf of weight $k$ on a presentation of a gerbe with atlas
$\mathcal{U} \to W$ is 
\begin{itemize}
\item
a sheaf $\mathcal{S}$ of $\mathcal{O}_{\mathcal{U}}$-modules over
$\mathcal{U}$
\item
an isomorphism $p_{0}^{*}\mathcal{S} \to 
(p_{1}^{*}\mathcal{S}) \otimes L^{k}$
\item
which satisfies the extra condition that the following diagram commutes

\begin{equation} \label{eqn:twistsheaf1}
\xymatrix@C-1pc{ p_{01}^{*}p_{0}^{*} \mathcal{S} \ar@<0ex>[r]^-{} & 
p_{01}^{*}(p_{1}^{*}\mathcal{S} \otimes L^{k}) \ar@<0ex>[r]^-{} &
(p_{01}^{*} p_{1}^{*}\mathcal{S}) \otimes p_{01}^{*}L^{k} \ar@<0ex>[r]^-{} & 
(p_{12}^{*}p_{0}^{*}\mathcal{S}) \otimes p_{01}^{*}L^{k} \ar@<0ex>[d]^-{}\\
p_{0}^{*}\mathcal{S} \ar@<0ex>[u]^-{} \ar@<0ex>[d]^-{} & & & 
p_{12}^{*}(p_{1}^{*}\mathcal{S} \otimes L^{k})\otimes p_{01}^{*}L^{k} 
\ar@<0ex>[d]^-{}
\\ p_{02}^{*}p_{0}^{*} \mathcal{S} \ar@<0ex>[r]^-{} 
& p_{02}^{*}(p_{1}^{*}\mathcal{S} \otimes L^{k}) \ar@<0ex>[r]^-{} 
& p_{2}^{*}\mathcal{S} \otimes p_{02}^{*}L^{k} & \ar@<0ex>[l]^-{} 
p_{2}^{*}\mathcal{S} \otimes p_{12}^{*}L^{k} \otimes p_{01}^{*}L^{k}}
\end{equation}
\end{itemize}

The sheaves on a presentation of a gerbe form an abelian category in 
the obvious way, a morphism between such sheaves $\mathcal{S}$ and 
$\mathcal{T}$ is simply a map $\mathcal{S} \to \mathcal{T}$ of 
sheaves on $\mathcal{U}$ such that the diagram below commutes
\begin{equation} \label{eqn:twistsheafmors}
\xymatrix{ p_{0}^{*}\mathcal{S} \ar@<0ex>[r]^-{} \ar@<0ex>[d]^-{} &
(p_{1}^{*}\mathcal{S}) \otimes L^{k} \ar@<0ex>[d]^-{} \\
p_{0}^{*}\mathcal{T} \ar@<0ex>[r]^-{}  &
(p_{1}^{*}\mathcal{T}) \otimes L^{k}}.
\end{equation}
In this way we get abelian categories $\text{Mod}({}_{(L, \theta,
\eta)} W, k)$ and $\text{Coh}({}_{(L, \theta, \eta)} W, k)$.  They can
also be seen as the categories of cartesian sheaves on the appropriate
semi-simplicial space.

For any presentation ${}_{(L, \theta, \eta)} W$ of a gerbe 
$\mathfrak{W} \to W$, we have isomorphisms
\[\text{Mod}({}_{(L, \theta, \eta)} W, k) \to \text{Mod}(\mathfrak{W}, k)
\]
and 
\[\text{Coh}({}_{(L, \theta, \eta)} W, k) \to \text{Coh}(\mathfrak{W}, k).
\]
\begin{remark}\label{rem:shv_on_se_pres}
Two presentations of gerbes on the same cover that are strongly
equivalent as defined in definition \ref{defi:strong_equiv_of_pres}
have equivalent categories of sheaves.
\end{remark}
\subsection{Derived Pushforward, Pullback, and Tensor Product} 
\label{ssec:DerivedFunctOnPresentations}
As mentioned before, there are some abstract ways to see that the
category of sheaves on a gerbe possess injective and flat objects.
To improve this slightly, we would like to see that there are flat and
injective objects in $\text{Mod}({}_{(L, \theta, \eta)} W) \cong
\text{Mod}(\mathfrak{W})$ which map to the same type of objects in
$\text{Mod}(\mathcal{U})$.  For this purpose consider the commutative
diagrams

\[
\xymatrix{
& \mathcal{U} \times_{\mathfrak{W}} \mathcal{U} \ar[dl]_-{\mathfrak{p}_{0}} 
\ar[dr]^-{\mathfrak{p}_{1}}
 & \\
\mathcal{U} \ar[dr]_-{\mathfrak{a}} & & \mathcal{U}
\ar[dl]^-{\mathfrak{a}} \\
& \mathfrak{W} &
}\hspace{15 mm} \text{and} \hspace{15 mm} 
\xymatrix{
& \mathcal{U} \times_{W} \mathcal{U} \ar[dl]_-{p_{0}} 
\ar[dr]^-{p_{1}}
 & \\
\mathcal{U} \ar[dr]_-{a} & & \mathcal{U}
\ar[dl]^-{a} \\
& W &
}.
\]
Notice that for any sheaf $\mathcal{M}$ of weight $k$ on
$\mathfrak{W}$ equation \ref{eqn:eqn-inj_inj} gives the existence of
an injection $\mathfrak{a}^{*}\mathcal{M} \to
\mathfrak{a}^{*}\mathfrak{a}_{*} \mathcal{I}$ for an injective sheaf
$\mathcal{I}$ on $\mathcal{U}$.  Flat base change then implies the
existence of an injection in $\text{Mod}({}_{(L, \theta, \eta)} W)$
\[\mathfrak{a}^{*}\mathcal{M} \to 
\mathfrak{a}^{*}\mathfrak{a}_{*}\mathcal{I} 
\cong \mathfrak{p}_{0*}\mathfrak{p}_{1}^{*}\mathcal{I} \cong 
p_{0*}(L^{k} \otimes p_{1}^{*} \mathcal{I})
\]  
Furthermore, $p_{0*}(L^{k} \otimes p_{1}^{*} \mathcal{I})$ is an injective
sheaf on $\mathcal{U}$.  Indeed since $p_{1}$ is a local isomorphism,
it is clear that $p_{1}^{*} \mathcal{I}$ is a locally injective
sheaf.  Then, by \cite[Proposition 2.4.10]{KS}, we conclude that
$p_{1}^{*} \mathcal{I}$ is injective.  For any sheaf $\mathcal{N}$,
we have $\text{Hom}(\mathcal{N}, L^{k} \otimes p_{1}^{*} \mathcal{I}) \cong
\text{Hom}(\mathcal{N} \otimes L^{-k}, p_{1}^{*} \mathcal{I})$ and
tensoring with $L^{-k}$ is exact, so $L^{k} \otimes p_{1}^{*}
\mathcal{I}$ is an injective ($\mathcal{O}_{\mathcal{U}^{1}}$-module).
Again, since the map $p_{0}: \mathcal{U}^{1} \to \mathcal{U}$ is a
local isomorphism we have that $\mathcal{O}_{\mathcal{U}^{1}} =
p_{0}^{-1}\mathcal{O}_{\mathcal{U}}$.  Finally, by
\cite[Proposition 2.4.1 (ii)]{KS} we have that $p_{0*}(L^{k} \otimes
p_{1}^{*} \mathcal{I})$ is injective.

Similarly, the equation \ref{eqn:eqn-flat_proj} gives a surjection
$\mathfrak{a}^{*} \mathfrak{a}_{!} \mathcal{F} \to
\mathfrak{a}^{*}\mathcal{M}$, for some flat sheaf $\mathcal{F}$ on 
$\mathcal{U}$.  Flat base change gives a surjection 
in $\text{Mod}({}_{(L, \theta, \eta)} W)$

\[p_{0!}(L^{k} \otimes p_{1}^{*}\mathcal{F}) 
\cong \mathfrak{p}_{0!}(\mathfrak{p}_{1}^{*}\mathcal{F}) \cong 
\mathfrak{a}^{*}\mathfrak{a}_{!} \mathcal{F} \to \mathfrak{a}^{*}\mathcal{M}.
\]
Now since $p_{1}$ is a local isomorphism, $p_{1}^{*} \mathcal{F}$ is
flat.  Then it is clear that $L^{k} \otimes p_{1}^{*} \mathcal{F}$ is also
flat, and flatness is preserved by $p_{0!}$, so we conclude that
$p_{0!}(L^{k} \otimes p_{1}^{*} \mathcal{F})$ is flat.

Suppose now that we have a commutative diagram

\[
\xymatrix{
\mathcal{A}  \ar@<0ex>[r]^-{a}
 \ar@<0ex>[d]^-{h} & M \ar@<0ex>[d] \\
\mathcal{B}  
\ar@<0ex>[r]^-{b} &
N
}
\]
where $\mathcal{A}$ is an atlas for $M$, $\mathcal{B}$ is an atlas for
$N$, and that the map $h$ induces a map $\tilde{h}: M \to N$.
Furthermore, suppose we are given a gerbe $\mathfrak{N}$ on $N$.
Chose a gerbe presentation ${}_{(L, \theta, \eta)} N$ of
$\mathfrak{N}$ using $\mathcal{B}$ and also consider the pullback
presentation ${}_{(h^{*}L, h^{*}\theta, h^{*}\eta)} N$ of
$h^{*}\mathfrak{N}$ on the atlas $\mathcal{A}$.  Then using
resolutions as we have described above (simultaneously flat or
injective on both the atlas and on the presentation), we see that the
following diagrams commute

\[
\xymatrix{ \sff{D}(\mathcal{A}) & \sff{D}({}_{(h^{*}L, h^{*}\theta,
h^{*}\eta)} M, k) \ar@<0ex>[l]^-{}  
& \ar@<0ex>[l]^-{\mathfrak{a}^{*}}  \sff{D}(\tilde{h}^{*}
\mathfrak{N}, k)  \\ \sff{D}(\mathcal{B}) \ar@<0ex>[u]^{Lh^{*}} & 
\sff{D}({}_{(L, \theta,\eta)} N, k) \ar@<0ex>[u]
\ar@<0ex>[l]^-{} & \ar@<0ex>[l]^-{\mathfrak{b}^{*}}  
\sff{D}(\mathfrak{N}, k) \ar@<0ex>[u]^{L \tilde{h}^{*}} }
\]
\medskip

\[
\xymatrix{ \sff{D}(\mathcal{A})  \ar@<0ex>[d]^{Rh_{*}} 
& \sff{D}({}_{(h^{*}L, h^{*}\theta,
h^{*}\eta)} M, k)  \ar@<0ex>[d] \ar@<0ex>[l]^-{}  
& \ar@<0ex>[l]^-{\mathfrak{a}^{*}}  \sff{D}(\tilde{h}^{*}
\mathfrak{N}, k) \ar@<0ex>[d]^{R \tilde{h}_{*}} \\ \sff{D}(\mathcal{B}) & 
\sff{D}({}_{(L, \theta,\eta)} N, k)
\ar@<0ex>[l]^-{} & \ar@<0ex>[l]^-{\mathfrak{b}^{*}}  
\sff{D}(\mathfrak{N}, k)  }.
\]

Also for any object $\mathcal{S}$ of $\sff{D}(\mathfrak{N}, c)$
and its associated object $\mathfrak{a}^{*}\mathcal{S}$ of
$\sff{D}({}_{(L, \theta, \eta)}N, c)$ we have the following
commutative diagram

\[
\xymatrix{ \sff{D}(\mathcal{B})  
\ar@<0ex>[d]^{LT_{\mathfrak{b}^{*}\mathcal{S}}} 
& \sff{D}({}_{(L, \theta,
\eta)} N, k)  \ar@<0ex>[d] \ar@<0ex>[l]^-{}  
& \ar@<0ex>[l]^-{\mathfrak{b}^{*}}  \sff{D}(
\mathfrak{N}, k) \ar@<0ex>[d]^{LT_{\mathcal{S}}} \\ \sff{D}(\mathcal{B}) & 
\sff{D}({}_{(L, \theta,\eta)} N, k +c)
\ar@<0ex>[l]^-{} & \ar@<0ex>[l]^-{\mathfrak{b}^{*}}  
\sff{D}(\mathfrak{N}, k+c)  }.
\]

Notice, that when thinking in terms of this presentation, we have for
every $k \in \mathbb{Z}$ an adjoint pair of functors 
$(\mathfrak{a}_{k}^{*},\mathfrak{a}_{k*})$, 
\begin{equation} \label{eqn:push2atlass}
\mathfrak{a}_{k}^{*}: \text{Mod}({}_{(L, \theta, \eta)}W,k) \to 
\text{Mod}(\mathcal{U}) 
\end{equation}
\begin{equation} \label{eqn:pull2atlass}
\mathfrak{a}_{k,*}: \text{Mod}(\mathcal{U}) \to 
\text{Mod}({}_{(L, \theta, \eta)}W,k)
\end{equation}
\noindent
Here $\mathfrak{a}_{k *}(\mathcal{S}) = p_{0*}(L^{k} \otimes
p_{1}^{*}\mathcal{S})$.  We will later use this in the form
\[
\mathfrak{a}_{-k *}(\mathcal{S}) = p_{1*}(L^{k} \otimes
p_{0}^{*}\mathcal{S}),
\]
when we will be interested in the case $k=1$.  Also note that
$\mathfrak{a}_{k}^{*}$ is the obvious forgetful map.

\section{Constructing New Presentations From Old} \label{ssec:FibTwist}
Let $\pi: X \to B$ be any map, $\mathcal{B}= \{ U_{i} \}$ an open
cover of $B$, with atlas 
\[U = \coprod_{i} U_{i} \to B,
\] and
$\mathcal{U}$ the pullback of this open cover by $\pi$.  In other
words we will use the atlas 
\[
a: \mathcal{U} = \coprod_{i} \pi^{-1}(U_{i}) \to X.
\]
\noindent
We will consider the space $X$ to be presented by the diagram

\[
\xymatrix{ 
\mathcal{U} \times_{X} 
\mathcal{U}\ar@<1ex>[r]^-{p^{1}_{0}}
\ar@<-2ex>[r]^-{p^{1}_{1}} & \mathcal{U} },
\]
where the two maps are the obvious projections
or in other words 

\[
\xymatrix{ 
\coprod_{ij} \pi^{-1}(U_{ij}) \ar@<1ex>[r]^-{p^{1}_{0}}
\ar@<-2ex>[r]^-{p^{1}_{1}} & \coprod_{i} \pi^{-1}(U_{i}) },
\]
where $p^{1}_{0}: \pi^{-1}(U_{ij}) \to \pi^{-1}(U_{i})$ and
$p^{1}_{1}: \pi^{-1}(U_{ij}) \to \pi^{-1}(U_{j})$ become the obvious
inclusion maps.  Let $\Delta : \mathcal{U} \to 
\mathcal{U} \times_{X} \mathcal{U}$
be the diagonal map.  This groupoid gives rise to a simplicial space
$\mathcal{U}^{\bullet}$ resolving $X$ over $B$ of the form

\[
\xymatrix{ 
\mathcal{U} \times_{X} 
\mathcal{U} \times_{X} \mathcal{U} \times_{X}
\mathcal{U}  \ar@<4ex>[r]^-{p^{3}_{012}} 
\ar@<1ex>[r]^-{p^{3}_{013}} \ar@<-2.5ex>[r]^-{p^{3}_{023}} 
\ar@<-5.5ex>[r]^-{p^{3}_{123}}
&
\mathcal{U} \times_{X} 
\mathcal{U} \times_{X} \mathcal{U} 
\ar@<2ex>[r]^-{p^{2}_{01}} \ar@<-1ex>[r]^-{p^{2}_{02}} 
\ar@<-4ex>[r]^-{p^{2}_{12}} & \mathcal{U} \times_{X} 
\mathcal{U} \ar@<1ex>[r]^-{p^{1}_{0}}
\ar@<-2ex>[r]^-{p^{1}_{1}} & \mathcal{U}}.
\]

Here we use the maps 
\[p^{n}_{i_{0}, \dots, i_{n}}: \mathcal{U}^{n+1} \to \mathcal{U}^{n}
\]
given by 
\[p^{n}_{i_{0}, \dots, i_{n-1}}(u_{0}, \dots, u_{n}) = 
(u_{i_{0}}, \dots, u_{i_{n-1}})
\]
We will not need the degeneracy maps too much, but they are the
obvious maps induced by the diagonal.  We will consider the sheaf of
groups $Aut(X/B)$ on $B$.  A section of this sheaf over an open set $U
\subset B$ is just an automorphism of $\pi^{-1}(U)$ which commutes
with the projection to the base.  To every element $f$ of
$\check{Z}^{1}_{\mathfrak{B}}(B, Aut(X/B))$ we will associate a new analytic
space $\pi_{f}: X_{f} \to B$.

Let $\pi^{n}:\mathcal{U}^{n} \to B$ be the obvious projection maps and
consider the automorphism $f$ of $\mathcal{U}^{1}=
\mathcal{U}\times_{X} \mathcal{U}$ satisfying $\pi^{1} = \pi^{1} \circ
f$.  We will denote the restriction of $f$ to the component
$\pi^{-1}(U_{ij})$ by $f_{ij}$.  On the space $\mathcal{U}^{a}$, we
write $f^{a}_{cd}: \mathcal{U}^{a} \to \mathcal{U}^{a}$, for the
restriction of $f$ to $\mathcal{U}^{a}$ via the map
\[
p^{a}_{cd}: \mathcal{U}^{a} \to \mathcal{U} \times_{X} \mathcal{U}.
\]
Therefore we have the equality $p_{cd}^{a} \circ f_{cd}^{a} = f
\circ p_{cd}^{a}$, and of course $f = f_{01}$.  The fact that $f$ is
in $\check{Z}^{1}(B, \mathcal{B}, Aut(X/B))$ just means that 
$f^{2}_{01} \circ f^{2}_{12} = f^{2}_{02}$ 
and that $f \circ \Delta = \Delta$.  We
will often just write this more conveniently as 
$f_{ij} \circ f_{jk} = f_{ik}$
on $\pi^{-1}(U_{ijk})$, and $f_{ii} = 1$ on $\pi^{-1}(U_{ii})
= \pi^{-1}(U_{i})$.

The map $f$ defines a new analytic space $X_{f}$ and the map
$\pi_{f}: X_{f} \to B$ is now defined by the presentation

\[\xymatrix{ 
\mathcal{U} \times_{X} \mathcal{U} \ar@<1.5ex>[r]^-{p^{1}_{0} \circ f}
\ar@<-2ex>[r]^-{p^{1}_{1} } & \mathcal{U} } \ \ \ \ \text{or} \ \ \ \
\xymatrix{ \coprod_{ij} \pi^{-1}(U_{ij}) \ar@<1.5ex>[r]^-{p^{1}_{0} \circ f}
\ar@<-2ex>[r]^-{p^{1}_{1}} & \coprod_{i} \pi^{-1}(U_{i}) }.
\]
In other words, $X_{f} = \mathcal{U}/\sim$ where we identify a point
$v \in \pi^{-1}(U_{j})$ with $f_{ij}(v) \in \pi^{-1}(U_{i})$ as long
as these points live over $U_{ij}$.

Notice that 
\[
\xymatrix{ 
\coprod_{ij} \pi^{-1}(U_{ij}) \ar@<1.5ex>[r]^-{p_{0}^{1} \circ f}
\ar@<-2ex>[r]^-{p^{1}_{1}} & \coprod_{i} \pi^{-1}(U_{i}) },
\]
is isomorphic to the diagram 
\[
\xymatrix{ 
\coprod_{ij} \pi_{f}^{-1}(U_{ij}) \ar@<1.5ex>[r]^-{}
\ar@<-2ex>[r]^-{} & \coprod_{i} \pi_{f}^{-1}(U_{i}) },
\]
(where we take the standard inclusions for the top and bottom arrows)
via the isomorphisms
\[\pi^{-1}(U_{ij}) \cong U_{i}
\times_{B} \pi^{-1}(U_{j}) = U_{i} 
\times_{B} \pi_{f}^{-1}(U_{j}) \cong \pi_{f}^{-1}(U_{ij}).
\]
This is perhaps more transparent in the other notation where we have
the commutative diagram
\[
\xymatrix{
\mathcal{U} \times_{X} \mathcal{U}  \ar@<1.5ex>[r]^-{p^{1}_{0} \circ f} 
\ar@<-2ex>[r]^-{p_{1}^{1}} \ar@<1ex>[d] & \mathcal{U} \\
\mathcal{U}_{f} \times_{X_{f}} \mathcal{U}_{f}   \ar@<1ex>[r]^-{}
\ar@<-1.5ex>[r]^-{} &
\mathcal{U}_{f} \ar@{=}[u]
}
\]
where the vertical arrow on the left hand side is the map $(p^{1}_{0}
\circ f, p_{1}^{1})$.  In what follows, we will usually use the top row
as a presentation of $X_{f}$.  Of course, an equivalence of 
two presentations of the twisted
fibrations $X_{f}$ and $X_{f'}$, given by $f$ and $f'$, is given by 
an automorphism $h$ of $\mathcal{U}$ which satisfies
$\pi \circ h = \pi$, and if we let $h_{c}$ denote the restriction of
$h$ to $\mathcal{U} \times_{X} \mathcal{U}$, via $p_{c}: \mathcal{U}
\times_{X} \mathcal{U} \to \mathcal{U}$,
\begin{equation} \label{eqn:twistequiv}
h_{0} \circ f = f' \circ h_{1}.
\end{equation}
This will often be written as $h_{i} \circ f_{ij} = f_{ij}' \circ h_{j}$ on 
the component $\pi^{-1}(U_{ij})$, where $h_{i}$ is just $h$ on the 
component $\pi^{-1}(U_{i})$.

Define 
\[m_{j}^{n}: (\mathcal{U}^{n})_{f} \to (\mathcal{U}^{n-1})_{f}
\] by
\[m^{n}_{0, \dots, n-1} = p_{0, \dots, n-1}^{n} \circ f_{n-1,n}.
\]
and for any $j$ with $0 \leq j \leq n-1$ $(i_{0}, \dots, i_{n-1})$
with $0 \leq i_{0} < i_{1} < \cdots i_{n-1} \leq n$, and $(i_{0},
\dots, i_{n-1}) \neq (0, \dots, n-1)$ let
\[m^{n}_{i_{0}, \dots, i_{n-1}} = p^{n}_{i_{0}, \dots, i_{n-1}}
\]
This defines a new simplicial manifold $(\mathcal{U}^{\bullet})_{f}$
whose first terms look like

\[
\xymatrix{ 
\mathcal{U} \times_{X} 
\mathcal{U} \times_{X} \mathcal{U} \times_{X}
\mathcal{U}  \ar@<4ex>[r]^-{p^{3}_{012} \circ f_{23}} 
\ar@<1ex>[r]^-{p^{3}_{013}} \ar@<-2.5ex>[r]^-{p^{3}_{023}} 
\ar@<-5.5ex>[r]^-{p^{3}_{123}}
&
\mathcal{U} \times_{X} 
\mathcal{U} \times_{X} \mathcal{U} 
\ar@<2.5ex>[r]^-{p^{2}_{01} \circ f_{12}} \ar@<-1ex>[r]^-{p^{2}_{02}} 
\ar@<-4ex>[r]^-{p^{2}_{12}} & \mathcal{U} \times_{X} 
\mathcal{U} \ar@<1ex>[r]^-{p^{1}_{0} \circ f}
\ar@<-2ex>[r]^-{p^{1}_{1}} & \mathcal{U}}
\]
and the degeneracy maps are induced by $f$ and the diagonal.  We have
an isomorphism of simplicial manifolds
\begin{equation} \label{eqn:isom_of_simp}
(\mathcal{U}^{\bullet})_{f} \cong (\mathcal{U}_{f})^{\bullet}.
\end{equation}  
Its $n$-th component 
\[(\mathcal{U}^{n})_{f} \to (\mathcal{U}_{f})^{n}
\]
is given by 

\[
\pi^{-1}(U_{t_{0}, t_{1}, \dots, t_{n}}) 
\cong U_{t_{0},\dots, t_{n-1}} \times_{B} \pi^{-1}(U_{t_{n}}) = 
U_{t_{0},\dots, t_{n-1}} \times_{B} \pi_{f}^{-1}(U_{t_{n}})
\cong \pi_{f}^{-1}(U_{t_{0}, t_{1}, \dots, t_{n}}).
\]

\subsection{Presentations of Gerbes on Twisted Fibrations} 
\label{ssec:gerbeontwist}
In order to define presentations of gerbes on the space $X_{f}$ we
will need to consider the terms
\[
\xymatrix{ 
\mathcal{U} \times_{X} 
\mathcal{U} \times_{X} \mathcal{U} \times_{X}
\mathcal{U}  \ar@<3.5ex>[r]^-{m_{012}} 
\ar@<1ex>[r]^-{m_{013}} \ar@<-2ex>[r]^-{m_{023}} 
\ar@<-4.5ex>[r]^-{m_{123}}
&
\mathcal{U} \times_{X} 
\mathcal{U} \times_{X} \mathcal{U} 
\ar@<2ex>[r]^-{m_{01}} \ar@<-.5ex>[r]^-{m_{02}} 
\ar@<-3ex>[r]^-{m_{12}} & \mathcal{U} \times_{X} 
\mathcal{U} \ar@<1ex>[r]^-{m_{0}}
\ar@<-1.5ex>[r]^-{m_{1}} & \mathcal{U}},
\]
coming from the resolution of $X_{f}$ over $B$ described above.
Due to the above remark, we see that a presentation of an
$\mathcal{O}^{\times}$ gerbe on $X_{f}$ is defined by a line bundle
$L$ on $\mathcal{U}\times_{X} \mathcal{U}$ together with isomorphisms
\[
\theta:m_{12}^{*}L \otimes m_{01}^{*}L \to m_{02}^{*}L
\]
and
\[
\eta: \mathcal{O}_{\mathcal{U}} \to \Delta^{*}L
\]
which satisfy the following conditions 

\begin{itemize}
\item \[
(m_{013}^{*} \theta) \circ (m_{123}^{*} \theta \otimes 1_{m_{01}^{*}L}) =
(m_{023}^{*} \theta) \circ (1_{m_{23}^{*}L} \otimes m_{012}^{*} \theta)
\]
as isomorphisms  
\[m_{23}^{*}L \otimes m_{12}^{*}L \otimes m_{01}^{*} L \to m_{03}^{*}L.
\]
Here we use 
\[m_{ij}: \mathcal{U} \times_{X} \mathcal{U} \times_{X}
\mathcal{U} \times_{X} \mathcal{U} \to \mathcal{U} \times_{X}
\mathcal{U}
\]
to denote the obvious maps.
\item \[(1_{\mathcal{U}} \times \Delta)^{*} 
\theta \circ (m_{1}^{*} \eta \times 1_{L}) = r
\]
\item \[(\Delta \times 1_{\mathcal{U}})^{*} \theta \circ 
(1_{L} \times m_{0}^{*} \eta) = l
\]
\end{itemize}
In other words, a presentation of a gerbe on a 
twisted version of $X$ is a quadruple 
\[
(f:\mathcal{U}^{1} \to \mathcal{U}^{1}, 
L \to \mathcal{U}^{1}, \theta, \eta)
\]
satisfying the above compatibilities.

We often write these as the isomorphisms
\[\theta_{ijk}: L_{jk} \otimes f_{jk}^{*}L_{ij} \to L_{ik}
\]
and 
\[\eta_{i}: \mathcal{O} \to L_{ii}
\]
which satisfy 
\begin{equation} \label{eqn:bigbad1}
\theta_{ijl} \circ (\theta_{jkl} \otimes 1_{L_{ij}})
= \theta_{ikl} \circ (1_{L_{kl}} \otimes f_{kl}^{*}\theta_{ijk}).
\end{equation}
\begin{equation} \label{eqn:bigbad2}
\theta_{ijj} \circ (\eta_{j} \otimes 1_{L_{ij}}) = l_{ij}
\end{equation}
and
\begin{equation} \label{eqn:bigbad3}
\theta_{iij} \circ(1_{L_{ij}} \otimes \eta_{j} ) = r_{ij}
\end{equation}
\begin{defi}\label{defi:equiv_pres_ger_tw_ver}
Two presentations ${}_{(L, \theta, \eta)}X_{f}$ and ${}_{(L', \theta',
\eta')}X_{f'}$ of gerbes on twisted versions of $X$ with respect to
the atlas $\mathcal{U} \to X$ are said to be strongly equivalent if
$f=f'$, and there is an isomorphism $L \to L'$ which intertwines
$\theta$ with $\theta'$ and $\eta$ with $\eta'$.
\end{defi}
Finally, let us comment on the notion of weak equivalence of two
presentations ${}_{(L, \theta, \eta)}X_{f}$ and ${}_{(L', \theta',
\eta')}X_{f'}$ of gerbes on twisted versions of $X$ over the same
cover.  Note that any equivalence $h$ of two twisted fibrations
$X_{f}$ and $X_{f'}$ gives rise to a map $\Xi(h)$ between
presentations of gerbes on $X_{f'}$ and presentations of gerbes on
$X_{f}$.  Guided by this observation we define a weak equivalence
between ${}_{(L, \theta, \eta)}X_{f}$ and ${}_{(L', \theta',
\eta')}X_{f'}$ to be a triple $(h; Q, \tau)$, where $h$ is an
equivalence between $f$ and $f'$, and the pair $(Q, \tau)$ is an
equivalence between the two gerbe presentations $(L, \theta, \eta)$
and $\Xi(h)(L', \theta', \eta')$ on $X_{f}$.  The map $\Xi(h)$ is
defined on objects by
\[(L', \theta', \eta') \mapsto 
(h_{1}^{*}L', h_{2}^{*}\theta', h_{0}^{*}\eta')
\]

\subsection{Classification of Gerbes and the 
Leray-Serre Spectral Sequence}\label{sec:classification}
Isomorphism classes of gerbes on the space $X_{f}$ are classified by
(the limit over covers $\mathfrak{C}$ of $X_{f}$) of the second total
cohomology group of the double complex
\[\check{C}_{\mathfrak{C}}^{q}((\mathcal{U}_{f})^{p}, \mathcal{O}^{\times})
\]
where $p \geq 0$ and $q \geq 0$.  We point out that the open sets in
$\mathfrak{C}$ are small and have nothing to do with $\mathcal{U}$ and
$\mathcal{U}_{f}$.  We will write this complex with $p$ increasing in
the horizontal direction and $q$ increasing in the vertical direction.
The rows only have cohomology in degree zero, so this double complex
is equivalent to the kernel column of the first horizontal map: the
column situated at $p = -1$ which is nothing but the \v{C}ech complex
$\check{C}_{\mathfrak{C}}^{\bullet}(X_{f}, \mathcal{O}^{\times})$ on
$X_{f}$, for the sheaf $\mathcal{O}^{\times}$.  In other words, we
have
\[H^{2}(X_{f}, \mathcal{O}^{\times}) \cong 
\mathbb{H}^{2}(\check{C}_{\mathfrak{C}}^{\bullet}
((\mathcal{U}_{f})^{\bullet}, \mathcal{O}^{\times})).
\]
The right hand should be thought of as follows: we can think of a
gerbe on $X_{f}$ as a gerbe on $\mathcal{U}_{f} =
(\mathcal{U}_{f})^{0}$ together with a twisted line bundle for the
difference of the pullback gerbes on $(\mathcal{U}_{f})^{1}$ (an
isomorphism of the two pullback gerbes) and
a relation between the pullback of these twisted line bundles to
$(\mathcal{U}_{f})^{2}$ which itself goes to the identity in
$(\mathcal{U}_{f})^{3}$.  In fact the double complex
$\check{C}_{\mathfrak{C}}^{q}((\mathcal{U}_{f})^{p},
\mathcal{O}^{\times})$ is the $E_{0}$ term of the Leray-Serre spectral
sequence which computes $H^{2}(X_{f}, \mathcal{O}^{\times})$ with
respect to the map $\pi_{f}: X_{f} \to B$.  We have
\[E_{f;0}^{p,q} = \check{C}_{\mathfrak{C}}^{q}
((\mathcal{U}_{f})^{p}, \mathcal{O}^{\times}) 
\]
\[E_{f;1}^{p,q} = H^{q}((\mathcal{U}_{f})^{p}, \mathcal{O}^{\times}) \cong 
\check{C}_{\mathfrak{B}}^{p}(B, R^{q}\pi_{f*} \mathcal{O}^{\times})
\]
and 
\[E_{f;2}^{p,q} = H^{p}(B, R^{q} \pi_{f*} \mathcal{O}^{\times}).
\]
In other words, we think of the Leray-Serre spectral sequence for
$\pi_{f}: X_{f} \to B$ as a spectral sequence associated with a
``\v{C}ech to \v{C}ech'' double complex where one of the covers is by
big opens and one is by small opens.

So, far everything we have said is general and applies
to any fibration $c: W \to B$, where we would replace
$\mathcal{U}_{f}^{\bullet}$ by the fiber products over $W$ of the disjoint
union of the pullback by $c$ of a cover of the base.

Now the isomorphism of simplicial manifolds
$(\mathcal{U}_{f})^{\bullet} \cong (\mathcal{U}^{\bullet})_{f}$ says
that we can classify gerbes on $X_{f}$ through data defined on $X$; we
have
\[H^{2}(X_{f}, \mathcal{O}^{\times}) \cong 
\mathbb{H}^{2}(\check{C}_{\mathfrak{C}}^{\bullet}
((\mathcal{U}_{f})^{\bullet}, \mathcal{O}^{\times})) \cong
\mathbb{H}^{2}(\check{C}_{\mathfrak{C}}^{\bullet}
((\mathcal{U}^{\bullet})_{f}, \mathcal{O}^{\times})). \] Of course we
will also use the isomorphism of simplicial manifolds
$(\mathcal{U}_{f})^{\bullet} \cong (\mathcal{U}^{\bullet})_{f}$ from
equation \ref{eqn:isom_of_simp}, to induce an isomorphism from the
spectral sequence $(E_{f;\bullet},d_{\bullet})$ to a new spectral
sequence $(E_{\bullet},d_{f;\bullet})$ where
\[E_{0}^{p,q} = 
\check{C}_{\mathfrak{C}}^{q}((\mathcal{U}^{p})_{f}, \mathcal{O}^{\times}) 
\]
\[E_{1}^{p,q} = H^{q}((\mathcal{U}^{p})_{f}, \mathcal{O}^{\times}) \cong 
\check{C}_{\mathfrak{B}}^{p}(B, R^{q}\pi_{f*} \mathcal{O}^{\times})
\]
and 
\[E_{2}^{p,q} = E_{f;2}^{p,q} = H^{p}(B, R^{q} \pi_{f*} \mathcal{O}^{\times}).
\]

\begin{rem} \label{rem-Classification_of_presentation}
Notice that the presentations of gerbes described in 
previous sections give classes in  
$\mathbb{H}^{2}(\check{C}^{\bullet}
((\mathcal{U}^{\bullet})_{f}, \mathcal{O}^{\times})) 
\cong H^{2}(X_{f}, \mathcal{O}^{\times})$ 
in the following manner.  
The line bundle $L$ is represented by an element of 
$\check{C}_{\mathfrak{C}}^{1}((\mathcal{U}^{1})_{f}, \mathcal{O}^{\times})$ 
and the the isomorphism $\theta$ is represented by an element of 
$\check{C}_{\mathfrak{C}}^{0}((\mathcal{U}^{2})_{f}, \mathcal{O}^{\times})$. 
\end{rem}

\section{Reinterpretation and Duality} \label{ssec:reinterpret}
A fibration $X \to B$ defines a stack of abelian categories on $B$, by
pushing forward the stack $\mathfrak{Mod}(X)$ from $X$.  To an open
set $U$ in $B$ we associate the stack of abelian categories
$\text{Mod}(\pi^{-1}(U))$.  On the other hand, given any gerbe over a
twisted version of $X$, which trivializes along the pullback of an
open cover of $B$, we can similarly associate a stack of abelian
categories on $B$.  Below, we explicitly re-express the data defining
a presentation of a gerbe over a twisted version of $X$ in terms of the
descent data of a stack of abelian categories on $B$, which is locally
on the base given as the stack on $B$ which was associated to $X$.  Of
course, this new stack defined by this descent data is isomorphic to the stack
of sheaves on the gerbe corresponding to the presentation.

 Consider the sheaf of groupoids given by $Pic(X/B)
\times Aut(X/B)$.  The multiplication rule for this sheaf is
\[
(L_{1}, f_{1})(L_{2}, f_{2}) = 
(L_{1} \otimes f_{1}^{*}L_{2}, f_{2} \circ f_{1})
\]

In these terms, a gerbe presentation on $X_{f}$ consists of isomorphisms 
\[(L_{jk}, f_{jk})(L_{ij}, f_{ij}) \to (L_{ik}, f_{ik})\]
or 
\[\theta_{ijk}: L_{jk} \otimes f_{jk}^{*}L_{ij} \to L_{ik}\]
and 
\[\eta_{i}:(\mathcal{O}, 1) \to (L_{ii}, 1)
\]
for which the diagram 

\[
\xymatrix{
(L_{kl}, f_{kl})(L_{jk}, f_{jk})(L_{ij}, f_{ij})  \ar@<0ex>[r]^-{} 
 \ar@<1ex>[d] & (L_{kl}, f_{kl})(L_{ik},f_{ik}) \ar@<1ex>[d] \\
(L_{jl}, f_{jl})(L_{ij}, f_{ij})  
\ar@<0ex>[r]^-{} &
(L_{il}, f_{il})
}
\]
is commutative.  Let 
\[A_{ij}: \text{Mod}(\pi^{-1}(U_{ij})) \to 
\text{Mod}(\pi^{-1}(U_{ij}))\]
be the functor given on objects by
\[A_{ij}(\mathcal{M}) = L_{ij} \otimes f_{ij}^{*}\mathcal{M}.
\]
We use 
$\btheta_{ijk}$ to denote the natural transformation 
\[
\btheta_{ijk}: A_{jk} \circ A_{ij}  \Rightarrow A_{ik}.
\]
Similarly, $\boldeta_{i}$ will denote the natural transformation 
\[\boldeta_{i}: 1_{ii}  \Rightarrow T_{L_{ii}} = T_{L_{ii}} \circ f_{ii}^{*} 
= A_{ii}.
\]
Thus we have reinterpreted a presentation of a gerbe as a collection
functors and natural transformations satisfying some commutative
diagrams.

\begin{lem} \label{lem:reenterpret}
A presentation of a gerbe over a twisted version of $X \to B$, on the
atlas $\mathcal{U} \to X$, defines a collection of functors $A_{ij}$
and natural transformations $\btheta_{ijk}$ and $\boldeta_{i}$ as
above, so that the following diagrams commute.

\begin{equation} \label{eqn:gerbe1}
\xymatrix{
A_{kl} \circ A_{jk} \circ A_{ij}  \ar@<0ex>[r]^-{A_{kl} \btheta_{ijk}} 
 \ar@<1ex>[d]^-{\btheta_{jkl} A_{ij}} & A_{kl} \circ A_{ik} 
\ar@<1ex>[d]^{\btheta_{ikl}} \\
A_{jl} \circ A_{ij}  
\ar@<0ex>[r]^-{\btheta_{ijl}} &
A_{il}
}
\end{equation}

\begin{equation} \label{eqn:gerbe2}
\xymatrix{1_{jj} \circ A_{ij} \ar@<0ex>[d]^-{\boldeta_{i} A_{ij}} 
\ar@<0ex>[r]^-{} & A_{ij} \\
A_{jj} \circ A_{ij} \ar@<-2ex>[ur]_-{\btheta_{ijj}} 
}
\end{equation}

\begin{equation} \label{eqn:gerbe3}
\xymatrix{
A_{ij} \circ 1_{ii} \ar@<0ex>[d]^-{1_{A_{ij}} \boldeta_{j}} \ar@<0ex>[r]^-{}
& A_{ij} 
\\
A_{ij} \circ A_{ii} \ar@<-2ex>[ur]_-{\btheta_{iij}} 
}
\end{equation}
\end{lem}

{\bf Proof}

The commutativity of these diagrams follows from equations
\ref{eqn:bigbad1},\ref{eqn:bigbad2} and \ref{eqn:bigbad3}.  Indeed we
need simply to observe that
\[
A_{kl} \btheta_{ijk} = 1_{L_{kl}} \otimes f_{kl}^{*} \theta_{ijk}
\] 
and 
\[\btheta_{jkl} A_{ij} = \theta_{jkl} \otimes 1_{L_{ij}}.
\] 
\ \hfill $\Box$
\begin{rem} \label{rem:onlyuptoisom}
Notice that if $(A, \btheta, \boldeta)$ are as above and $A' \in
Aut(\text{Coh}(\mathcal{U}^{1}))$ is a functor that is isomorphic to
$A$ in a way compatible with restriction on base, then there exist
canonical natural transformations $\btheta'$, and $\boldeta'$ such
that $(A', \btheta', \boldeta')$ satisfy the same diagrams as do $(A,
\btheta, \boldeta)$.  Therefore, we could have taken $A$ to be any
automorphism of stacks isomorphic to $T_{L} \circ f^{*}$.
\end{rem}

\begin{defi} \label{def:strongequiv}
We call two triples $(A, \btheta, \boldeta)$, and $(A',\btheta',
\boldeta')$ strongly equivalent if there is a natural transformation
between $A$ and $A'$ which is compatible with restriction along the
base, which intertwines $\btheta$ with $\btheta'$, and also
intertwines $\boldeta$ with $\boldeta'$.
\end{defi}

A weak equivalence of triples $(A, \btheta, \boldeta)$,
and $(A',\btheta', \boldeta')$ is given by a pair $(B, \bupsilon)$,
where $B:\text{Mod}(\mathcal{U}) \to \text{Mod}(\mathcal{U})$ is an
autoequivalence and
\begin{equation} \label{eqn:upsilon}
\bupsilon_{ij}: A_{ij} \circ B_{i} \Rightarrow B_{j} \circ A'_{ij} 
\end{equation}
is an invertible natural transformation 
which satisfies the equality of natural transformations 
\[\bupsilon_{ik} \circ (\btheta_{ijk}B_{i}) = (B_{k} \btheta_{ijk}') 
\circ (\bupsilon_{jk} A_{ij}') \circ (A_{jk} \bupsilon_{ij})
\] 
as natural transformations 
\[A_{jk} \circ A_{ij} \circ B_{i} \Rightarrow B_{k} \circ A_{ik}'.
\]
Indeed, we can recognize the equation \ref{eqn:upsilon} as containing
the maps
\[
L_{ij} \otimes f_{ij}^{*}Q_{i} \to Q_{j} \otimes h_{j}^{*} L'_{ij}.
\]
\begin{rem} \label{rem:fullyfaithfull}
The lemma \ref{lem:reenterpret} gives a fully faithful functor from
the category whose objects are presentations ${}_{(L, \theta, \eta)}
X_{f}$ of gerbes on twisted versions of $X \to B$ using the atlas
$\mathcal{U} \to X$ and whose morphisms are weak/strong equivalences,
to the category whose objects are triples $(A, \btheta, \boldeta)$ and
whose morphisms are weak/strong equivalences.
\end{rem}
\begin{rem}\label{rem:reints}

Any triple $(A, \btheta, \boldeta)$ defines an abelian category whose
objects are sheaves $\mathcal{S}$ on $\mathcal{U}$, along with
isomorphisms
\[A(p_{0}^{*}\mathcal{S}) \to  A(p_{1}^{*}\mathcal{S}) 
\hspace{15 mm} \text{or}  \hspace{15 mm} 
A_{ij}(\mathcal{S}_{i}) \to \mathcal{S}_{j}
\]
such that 

\begin{equation} \label{eqn:reenterprettwist}
\xymatrix{
(p_{12}^{*}A \circ p_{01}^{*}A)p_{0}^{*}\mathcal{S} \ar@<0ex>[r]^-{} 
\ar@<0ex>[d]^-{} & (p_{12}^{*}A)(p_{2}^{*}\mathcal{S}) \ar@<0ex>[d]^-{}
\\
(p_{02}^{*}A)(p_{0}^{*}\mathcal{S}) \ar@<0ex>[r]^-{} & p_{2}^{*}\mathcal{S}
}
\hspace{15 mm} \text{or}  \hspace{15 mm}
\xymatrix{
(A_{jk} \circ A_{ij})\mathcal{S}_{i} \ar@<0ex>[r]^-{} 
\ar@<0ex>[d]^-{} & A_{jk}(\mathcal{S}_{j}) \ar@<0ex>[d]^-{}
\\
A_{ik}(\mathcal{S}_{i}) \ar@<0ex>[r]^-{} & \mathcal{S}_{k}
}
\end{equation}
commutes, and morphisms are maps $\mathcal{S} \to \mathcal{T}$ such
that the diagram
\[
\xymatrix{
A(p_{0}^{*} \mathcal{S}) \ar@<0ex>[r]^-{} 
\ar@<0ex>[d]^-{} & p_{1}^{*} \mathcal{S} \ar@<0ex>[d]^-{}
\\
A(p_{0}^{*} \mathcal{T}) \ar@<0ex>[r]^-{} & p_{1}^{*} \mathcal{T}
}\hspace{15 mm} \text{or}  \hspace{15 mm}
\xymatrix{
A_{ij}(\mathcal{S}_{i}) \ar@<0ex>[r]^-{} 
\ar@<0ex>[d]^-{} & \mathcal{S}_{j} \ar@<0ex>[d]^-{}
\\
A_{ij}(\mathcal{T}_{i}) \ar@<0ex>[r]^-{} & \mathcal{T}_{j}
}
\]
commutes.  Any strong equivalence between triples $(A, \btheta,
\boldeta)$ and $(A', \btheta', \boldeta')$ induces an isomorphism of
the associated abelian categories.  If the triple $(A, \btheta,
\boldeta)$ comes from a gerbe presentation ${}_{(L, \theta, \eta)}
X_{f}$ then the category we have described is isomorphic to the
category of sheaves of weight $(-1)$ on the presentation ${}_{(L,
\theta, \eta)} X_{f}$.
In the future, we will make this identification implicit.
\end{rem}
\begin{defi} \label{defi:compatable}

Suppose we are given maps $\pi: X \to B$ and $\rho: Y \to B$, and an
object $\mathcal{P}$ in $\sff{D}(X \times_{B} Y)$ giving a derived
equivalence $\Phi: \sff{D}(Y) \to \sff{D}(X)$.  We call a gerbe
presentation ${}_{(L, \theta, \eta)} X_{f}$ $\Phi$-compatible if we
have
\[A \circ \Phi = \Phi \circ C 
\]
where $A \in Aut(\text{Coh}(\mathcal{U}^{1}))$ corresponds to ${}_{(M,
\theta, \eta)} X_{f}$ and $C\in Aut(\text{Coh}(\mathcal{V}^{1}))$, is
isomorphic to the functor $\mathcal{S} \to M \otimes
g^{*}\mathcal{S}$, for some line bundle $M$, and automorphism $g$
preserving the projection to the base.  
\\
\noindent
\\ We call a gerbe $\mathfrak{X} \to B$ over a twisted version of $X$
$\Phi$-compatible if it admits a $\Phi$-compatible presentation
${}_{(L, \theta, \eta)} X_{f}$.
\end{defi}

\subsection{The Dual Gerbe on the Dual Fibration}
Using our description of presentations of gerbes on twisted
versions of $X$, we can now employ the base changed derived equivalence
\[\Phi_{U^{1}} :\sff{D}(\mathcal{V}) \to \sff{D}(\mathcal{U})
\]
to map $\Phi$-compatible presentations of gerbes on twisted versions
of $\pi: X \to B$ to $\Psi$-compatible presentations of gerbes on
twisted versions of $\rho: Y \to B$, where we have chosen an inverse
functor $\Psi: \sff{D}(Y) \to \sff{D}(X)$ to $\Phi$.  It may help to
keep in mind the commutative diagram

\[
\xymatrix{
& \mathcal{U} \times_{U} \mathcal{V}  \ar[dl]_-{\tilde{\rho}_{U}}  
\ar[dr]^-{\tilde{\pi}_{U}}
\ar[dd]_-{r}& \\
\mathcal{U} \ar[dd]_-{a} 
& & \mathcal{V} \ar[dd]_-{b} \\
& X \times_{B} Y \ar[dl]_-{\tilde{\rho}} \ar[dr]^-{\tilde{\pi}}
 & \\
X \ar[dr]_-{\pi} & & Y
\ar[dl]^-{\rho} \\
& B &
}
\]
\noindent
on which our presentations are based.

Consider the functor 
\[
\Upsilon: Aut (\text{Coh}(\mathcal{U}^{1})) \to Aut(\sff{D}(\mathcal{V}^{1}))
\]
given by 
\[\Upsilon: A \mapsto \Psi_{U^{1}} \circ A \circ \Phi_{U^{1}}.
\]
We will denote this on elements $A_{ij}$ by the equation 
\begin{equation} \label{eqn:twistsheaf1}
\Upsilon(A_{ij}) = \Psi \circ A_{ij} \circ \Phi.
\end{equation}
The functor $\Upsilon(A_{ij})$ satisfies  
\[\Upsilon(A_{ij}) \cong T_{M_{ij}} \circ g_{ij}^{*},\] for 
some line bundles $M_{ij}$ on $\rho^{-1}(U_{ij})$, and some fiber preserving
automorphisms $g_{ij}$ of $\rho^{-1}(U_{ij})$.  Here the isomorphism
respects the pullback along open sets in the base.  In other words,
this is an isomorphism of automorphisms of the pushforward of the
stack of coherent sheaves on $Y$ to the base $B$.

Let $Aut(\text{Mod}(\mathcal{U}^{1}))^{\Phi}$ be the sub group-groupoid of
$Aut(\text{Mod}(\mathcal{U}^{1}))$ consisting of the elements $A$ for which
$\Upsilon(A)$ is isomorphic to $T_{M} \circ g^{*}$ for some line
bundle $M$ and projection preserving automorphism $g$.  We get an
induced isomorphism of categories, where the morphisms are strong
equivalences
\[
\Upsilon: Aut(\text{Mod}(\mathcal{U}^{1}))^{\Phi} \to  
Aut(\text{Mod}(\mathcal{V}^{1}))^{\Psi}.
\]
By applying $\Upsilon$ to the diagram \ref{eqn:gerbe1} we obtain the 
following commutative diagram.
\begin{equation} \label{eqn:upsilongerbe}
\xymatrix
{
\Upsilon(A_{kl} \circ A_{jk} \circ A_{ij})  
\ar@<0ex>[r]^-{\Upsilon(A_{kl} \btheta_{ijk})} 
 \ar@<1ex>[d]^-{\Upsilon(\btheta_{jkl} A_{ij})} & 
\Upsilon (A_{kl} \circ A_{ik}) 
\ar@<1ex>[d]^{\Upsilon(\btheta_{ikl})} \\
\Upsilon(A_{jl} \circ A_{ij})  
\ar@<0ex>[r]^-{\Upsilon(\btheta_{ijl})} &
\Upsilon(A_{il})
}
\end{equation}

This is not quite what we want, but if we attach to this 
diagram three other commuting squares, we get the outer diagram below, 
which is the one that we want.   
\begin{equation} \label{eqn:bigupsilongerbe}
\xymatrix{
\Upsilon(A_{kl}) \circ \Upsilon(A_{jk}) \circ \Upsilon(A_{ij}) 
\ar@<0ex>[r] \ar@<1ex>[d]
& \Upsilon(A_{kl}) \circ \Upsilon(A_{jk} \circ A_{ij}) 
\ar@<0ex>[r] \ar@<1ex>[d] &
\Upsilon(A_{kl}) \circ \Upsilon(A_{ik}) \ar@<0ex>[d]^-{} \\
\Upsilon(A_{kl} \circ A_{jk}) \circ \Upsilon(A_{ij}) 
\ar@<0ex>[r] \ar@<1ex>[d]
& \Upsilon(A_{kl} \circ A_{jk} \circ A_{ij})  
\ar@<0ex>[r]^-{\Upsilon(A_{kl} \btheta_{ijk})} 
\ar@<1ex>[d]^-{\Upsilon(\btheta_{jkl} A_{ij})} & 
\Upsilon (A_{kl} \circ A_{ik}) 
\ar@<1ex>[d]^{\Upsilon(\btheta_{ikl})} \\
\Upsilon(A_{jl}) \circ \Upsilon(A_{ij}) \ar@<0ex>[r] 
& \Upsilon(A_{jl} \circ A_{ij})  
\ar@<0ex>[r]^-{\Upsilon(\btheta_{ijl})} &
\Upsilon(A_{il})
}
\end{equation}
We now rewrite this, defining $\bT_{ikl}$ to be the right hand 
side of the above diagram: 
\begin{equation} \label{eqn:Newgerbe1}
\xymatrix
{
\Upsilon(A_{kl}) \circ \Upsilon(A_{jk}) \circ \Upsilon(A_{ij})
  \ar@<0ex>[r]^-{\Upsilon(A_{kl}) \bT_{ijk}} 
 \ar@<1ex>[d]_-{\bT_{jkl} \Upsilon(A_{ij})} & \Upsilon(A_{kl}) 
\circ \Upsilon(A_{ik}) 
\ar@<1ex>[d]^{\bT_{ikl}} \\
\Upsilon(A_{jl}) \circ \Upsilon(A_{ij})  
\ar@<0ex>[r]^-{\bT_{ijl}} &
\Upsilon(A_{il}).
}
\end{equation}

Similarly by applying $\Upsilon$ to the diagrams $\ref{eqn:gerbe2}$ 
and $\ref{eqn:gerbe3}$ we get the diagrams

\begin{equation} \label{eqn:Upsilongerbe2}
\xymatrix{\Upsilon(1_{jj} \circ A_{ij}) 
\ar@<0ex>[d]_-{\Upsilon(\boldeta_{i} A_{ij})} 
\ar@<0ex>[r]^-{} & \Upsilon(A_{ij}) \\
\Upsilon(A_{jj} \circ A_{ij}) \ar@<-1ex>[ur]_-{\Upsilon(\btheta_{jji})} 
}
\end{equation}
and

\begin{equation} \label{eqn:Upsilongerbe3}
\xymatrix{
\Upsilon(A_{ij} \circ 1_{ii}) 
\ar@<0ex>[d]_-{\Upsilon(1_{A_{ij}} \boldeta_{j})} \ar@<0ex>[r]^-{}
& \Upsilon(A_{ij}) 
\\
\Upsilon(A_{ij} \circ A_{ii}) \ar@<-1ex>[ur]_-{\Upsilon(\btheta_{jji})} 
}
\end{equation}
By applying $\Upsilon$ to the natural transformation $\boldeta_{i}: 1_{ii}
\Rightarrow A_{ii}$, we can define $\bzeta_{i}$ as the composition 
$1_{ii} \Rightarrow \Upsilon(1_{ii}) \Rightarrow \Upsilon(A_{ii})$.
Using $\bzeta_{i}: 1_{ii} \Rightarrow \Upsilon(A_{ii})$ we can complete the
above diagram as
\begin{equation} \label{eqn:Upsilongerbe2complete}
\xymatrix{
1_{jj} \circ \Upsilon(A_{ij}) 
\ar@<0ex>[r]^-{} \ar@<0ex>[d]^-{} & \Upsilon(1_{jj} \circ A_{ij}) 
\ar@<0ex>[d]_-{\Upsilon(\boldeta_{j} A_{ij})} 
\ar@<0ex>[r]^-{} & \Upsilon(A_{ij}) \\
\Upsilon(A_{jj}) \circ \Upsilon(A_{ij})\ar@<0ex>[r]^-{} & 
\Upsilon(A_{jj} \circ A_{ij}) \ar@<-1ex>[ur]_-{\Upsilon(\btheta_{ijj})} 
}
\end{equation}

and

\begin{equation} \label{eqn:blah}
\xymatrix{
\Upsilon(A_{ij}) \circ 1_{ii}
\ar@<0ex>[r]^-{} \ar@<0ex>[d]^-{} & \Upsilon(A_{ij} \circ 1_{ii}) 
\ar@<0ex>[d]_-{\Upsilon(A_{ij} \boldeta_{i})} 
\ar@<0ex>[r]^-{} & \Upsilon(A_{ij}) \\
 \Upsilon(A_{ij}) \circ \Upsilon(A_{ii}) \ar@<0ex>[r]^-{} & 
\Upsilon(A_{ij} \circ A_{ii}) \ar@<-1ex>[ur]_-{\Upsilon(\btheta_{iij})} 
}
\end{equation}

In other words, we obtain the commutative diagrams

\begin{equation} \label{eqn:Newgerbe2}
\xymatrix{
1_{jj} \circ \Upsilon(A_{ij})
\ar@<0ex>[r]^-{} \ar@<0ex>[d]_-{\bzeta_{j}\Upsilon(A_{ij})} 
& \Upsilon(A_{ij}) \\
\Upsilon(A_{jj}) \circ \Upsilon(A_{ij}) 
\ar@<0ex>[ur]_-{\bT_{ijj}} 
}
\end{equation}

and

\begin{equation} \label{eqn:Newgerbe3}
\xymatrix{
\Upsilon(A_{ij}) \circ 1_{ii}
\ar@<0ex>[r]^-{} \ar@<0ex>[d]_-{\Upsilon(A_{ij}) \bzeta_{i}} 
& \Upsilon(A_{ij}) \\
\Upsilon(A_{ij}) \circ \Upsilon(A_{ii}) 
\ar@<0ex>[ur]_-{\bT_{iij}} 
}
\end{equation}

Notice that we end up with the same type of diagrams
$\ref{eqn:Newgerbe1}, \ref{eqn:Newgerbe2}, \ref{eqn:Newgerbe3}$ that
we started with on the other side $\ref{eqn:gerbe1}, \ref{eqn:gerbe2},
\ref{eqn:gerbe3}$ .  Finally, we apply the diagrams
$\ref{eqn:Newgerbe1}, \ref{eqn:Newgerbe2}, \ref{eqn:Newgerbe3}$ to the
structure sheaves of points $\mathbb{C}_{y}$ on $Y$, and to the
structure sheaf $\mathcal{O}_{Y}$ of $Y$.  The application to points
shows that $g$ is a cocycle.  Application to the structure sheaf gives
a presentation of a gerbe on $Y_{g}$.  This gives rise to a map
\[
b:\mathcal{V} = \coprod_{i} \rho^{-1}(U_{i}) \to Y,
\]
an automorphism $g: \mathcal{V}^{1} \to \mathcal{V}^{1}$, a line bundle
$M \to \mathcal{V}^{1}$, and isomorphisms
\[T_{ijk}: M_{jk} \otimes g_{jk}^{*}M_{ij} \to M_{ik}
\]
and
\[\zeta_{i} : \mathcal{O} \to M_{ii}
\]
which satisfy

\begin{equation} \label{eqn:newbigbad1}
T_{ijl} \circ (T_{jkl} \otimes 1_{M_{ij}})
= T_{ikl} \circ (1_{M_{kl}} \otimes g_{kl}^{*}T_{ijk}).
\end{equation}
\begin{equation} \label{eqn:newbigbad2}
T_{ijj} \circ (\zeta_{j} \otimes 1_{M_{ij}}) = l_{ij}
\end{equation}
and
\begin{equation} \label{eqn:newbigbad3}
T_{iij} \circ(1_{M_{ij}} \otimes \zeta_{j} ) = r_{ij}.
\end{equation} 

Thus, we have found a presentation $(g :
\mathcal{V}^{1} \to \mathcal{V}^{1}, 
M \to \mathcal{V}^{1}, T, \zeta)$ of a gerbe on
a twisted version of $Y$.  This gerbe presentation ${}_{(M,T, \zeta)}
Y_{g} \to Y_{g}$ over a twisted version of $Y$, was found starting from the
presentation of the gerbe ${}_{(L,\theta, \eta)} X_{f} \to X_{f}$ over
a twisted version of $X$.  We state this as the following
\begin{lem} \label{lem:thedual}
Given maps $\pi:X \to B$ and $\rho: Y \to B$ and an object
$\mathcal{P} \in \sff{D}(X \times_{B} Y)$ giving a derived equivalence
$\Phi: \sff{D}(Y) \to \sff{D}(X)$ and a quasi-inverse $\Psi$, there is
a correspondence between $\Phi$-compatible presentations of gerbes on
twisted versions of $X$ and $\Psi$-compatible presentations of gerbes
on twisted versions of $Y$.  Using strong equivalences as morphisms,
this gives an isomorphism of categories.
\end{lem}
{\bf Proof} Above.\ \hfill $\Box$
\begin{defi} \label{defi:dual}
Suppose $\mathfrak{X}$ is a $\Phi$-compatible gerbe on a twisted
version of $X$ and $\mathfrak{Y}$ is the corresponding gerbe on a
twisted version of $Y$ produced as described above.  Then we say that 
$\mathfrak{Y}$ is $\Phi$-dual to $\mathfrak{X}$.  $\mathfrak{Y}$ is
defined by a $\Phi$-dual presentation ${}_{(M,T,\zeta)} Y_{g}$ to a
$\Phi$-compatible presentation ${}_{(L,\theta,\eta)} X_{f}$ of
$\mathfrak{X}$.
\end{defi}

Notice that the definition of $T$ is encoded in the diagram
\begin{equation} \label{eqn:twistsheaf2}
\xymatrix{
\Psi \circ A_{jk} \circ \Phi \circ \Psi \circ A_{ij} \circ \Phi
 \ar@<0ex>[dr]^-{\bT_{ijk}} &  \\
\Psi \circ A_{jk} \circ A_{ij} \circ \Phi  \ar@<0ex>[u]^-{}
\ar@<0ex>[r]^-{} &
\Psi \circ A_{ik} \circ \Phi.
}
\end{equation}

\subsection{The New Functor}

It turns out that equations \ref{eqn:twistsheaf1} and
\ref{eqn:twistsheaf2} secretly encode the description of a new functor. 
We have the following lemma.
\begin{lem} \label{lem:newfunctor}
Suppose we are given compact, connected, complex manifolds $X$ and
$Y$, maps $\pi:X \to B$ and $\rho: Y \to B$, where $\rho$ is flat and
an object \[
\mathcal{P} \in \text{Coh}(X \times_{B} Y),
\] flat over $Y$
giving a derived equivalence $\Phi: \sff{D}^{b}_{c}(Y) \to
\sff{D}^{b}_{c}(X)$ and a quasi-inverse equivalence $\Psi$.  Then for
every $\Phi$-compatible gerbe presentation ${}_{(L, \theta, \eta)}
X_{f}$ and $\Phi$-dual gerbe presentation ${}_{(M, T, \zeta)} Y_{g}$,
we can define a functor
\[\widetilde{\Phi} :\sff{D}^{b}_{c}({}_{(M, T, \zeta)} Y_{g}, -1) 
\to \sff{D}^{b}_{c}({}_{(L, \theta, \eta)} X_{f},-1) \]
\end{lem}

{\bf Proof}

In order to see this, we must first
manipulate the equation \ref{eqn:twistsheaf2} a little.  By applying
$\Phi$ on the left to the diagram \ref{eqn:twistsheaf2} and attaching a
commutative lower square, we get
\begin{equation} \label{eqn:mhotwistsheaf2}
\xymatrix{
\Phi \circ \Psi \circ A_{jk} \circ \Phi 
\circ \Psi \circ A_{ij} \circ \Phi 
 \ar@<0ex>[dr]^-{\Phi \bT_{ijk}} &  \\
\Phi \circ \Psi \circ A_{jk} \circ A_{ij} \circ \Phi  
\ar@<0ex>[r]^-{} \ar@<0ex>[u]^-{} &
\Phi \circ \Psi \circ A_{ik} \circ \Phi  & \\
A_{jk} \circ A_{ij} \circ \Phi \ar@<0ex>[r]^-{\btheta_{ijk} \Phi} 
\ar@<0ex>[u]^-{}
& A_{ik} \circ \Phi \ar@<0ex>[u]^-{} }.
\end{equation}
Taking the outer part of this we arrive at the commutative diagram 
\begin{equation} \label{eqn:anotherthing}
\xymatrix{
\Phi \circ \Upsilon(A_{jk}) \circ \Upsilon(A_{ij})
\ar@<1ex>[d]^-{\Phi \bT_{ijk}} & \ar@<0ex>[l]^-{}  A_{jk} 
\circ A_{ij} \circ \Phi
\ar@<1ex>[d]^{\btheta_{ijk} \Phi} \\
\Phi \circ \Upsilon(A_{ik})  
& \ar@<0ex>[l]^-{}  A_{ik} \circ \Phi
}
\end{equation}
\begin{rem}\label{rem:alass}
Every corner of this square is a functor implemented by a sheaf on the
appropriate product.  If we could find unique sheaf isomorphisms
implementing the natural equivalences indicated by the arrows, this
diagram would describe a sheaf of weight $(-1,1)$ on the product of
our gerbe presentations and it would be easy to show that such a sheaf
implements an equivalence using the convolution product.  However, it is not
clear to the author how to produce these isomorphisms, and so we must
follow a different path.
\end{rem}

We start by defining a functor 
\[\widetilde{\phi}:
\text{Mod}({}_{(M,T,\zeta)} Y_{g},
-1) \to \text{Mod}({}_{(L,\theta,\eta)} X_{f}, -1).
\]
We will sometimes write this functor as 
\[\widetilde{\phi}(\mathcal{S}, \kappa: \Upsilon(A)p_{0}^{*}\mathcal{S} 
\to p_{1}^{*} \mathcal{S}) = 
(\widetilde{\phi}(\mathcal{S}),\widetilde{\phi}(\kappa)).
\]
Here we define $\widetilde{\phi}(\mathcal{S}) = 
\mathcal{H}^{0}(\Phi_{U}(\mathcal{S}))$ and 
$\widetilde{\phi}(\kappa):A(p_{0}^{*}\widetilde{\phi}(\mathcal{S})) 
\to p_{1}^{*}\widetilde{\phi}(\mathcal{S})$ is defined as the composition
\[A(p_{0}^{*}(\mathcal{H}^{0}(\Phi_{U}(\mathcal{S})))) 
\to \mathcal{H}^{0}(Ap_{0}^{*}(\Phi_{U}(\mathcal{S}))) 
\to \mathcal{H}^{0}(A(\Phi_{U^{1}} p_{0}^{*}(\mathcal{S}))) 
\to \mathcal{H}^{0}(\Phi_{U^{1}}(\Upsilon(A)p_{0}^{*}\mathcal{S})) 
\to 
\] 
\[\to 
\mathcal{H}^{0}(\Phi_{U^{1}}(p_{1}^{*}\mathcal{S})) 
\to p_{1}^{*}(\mathcal{H}^{0}(\Phi_{U} \mathcal{S})).
\]
For the sake of sanity, we will change notation, skip the second and
second to last term.  We write the object of
$\text{Mod}({}_{(M,T,\zeta)} Y_{g}, -1)$ as 
\[(\{ \mathcal{S}_{i} \},
\{\kappa_{ij} : \Upsilon(A_{ij}) \mathcal{S}_{i} \to \mathcal{S}_{j}
\}).
\]
It is sent by
$\widetilde{\phi}$ to 
\[(\{\mathcal{H}^{0}(\Phi \mathcal{S}_{i})\},
\{\lambda_{ij} \}),
\] where $\lambda_{ij}$ is defined by the
composition
\[A_{ij} \mathcal{H}^{0}(\Phi \mathcal{S}_{i}) \to 
\mathcal{H}^{0} (A_{ij}\Phi \mathcal{S}_{i}) \to 
\mathcal{H}^{0}(\Phi \Upsilon(A_{ij}) \mathcal{S}_{i}) 
\to \mathcal{H}^{0}(\Phi \mathcal{S}_{j})
\]
\[\widetilde{\phi}(\{ \mathcal{S}_{i} \}, \{\kappa_{ij}\}) = 
(\{\mathcal{H}^{0}(\Phi \mathcal{S}_{i})\},  \{\lambda_{ij} \}). 
\]

The fact that this is an object of 
$\text{Mod}({}_{(M,T,\zeta)} Y_{g}, -1)$ rests on the commutativity of
the outer square in the following:
\[
\xymatrix@C-1pc{A_{jk} A_{ij} \mathcal{H}^{0}(\Phi \mathcal{S}_{i}) 
\ar[r]^-{} \ar[ddd]^-{}&
A_{jk} \mathcal{H}^{0}(A_{ij} \Phi \mathcal{S}_{i}) 
\ar[r]^-{} \ar[dd]^-{}&
A_{jk} \mathcal{H}^{0}(\Phi \Upsilon(A_{ij})\mathcal{S}_{i}) 
\ar[r]^-{} \ar[d]^-{} &
A_{jk} \mathcal{H}^{0}(\Phi \mathcal{S}_{j}) \ar[d]^-{} \\
& & \mathcal{H}^{0}(A_{jk} \Phi \Upsilon(A_{ij})\mathcal{S}_{i}) 
\ar[r]^-{} \ar[d]^-{}
& \mathcal{H}^{0}(A_{jk} \Phi \mathcal{S}_{j}) \ar[d]^-{} \\
& \mathcal{H}^{0}(A_{jk}  A_{ij} \Phi \mathcal{S}_{i}) 
\ar[r]^-{} \ar[d]^-{} & 
\mathcal{H}^{0}(\Phi \Upsilon(A_{jk}) \Upsilon(A_{ij}) \mathcal{S}_{i}) 
\ar[r]^-{} \ar[d]^-{}& 
\mathcal{H}^{0}(\Phi \Upsilon(A_{jk}) \mathcal{S}_{j}) \ar[d]^-{} \\
A_{ik} \mathcal{H}^{0}(\Phi \mathcal{S}_{i})\ar[r]^-{} & 
\mathcal{H}^{0}(A_{ik} \Phi \mathcal{S}_{i}) \ar[r]^-{}&
\mathcal{H}^{0}(\Phi \Upsilon(A_{ik}) \mathcal{S}_{i}) \ar[r]^-{}&
\mathcal{H}^{0}(\Phi \mathcal{S}_{k}).
}
\]
\noindent
We can see that the outer square is commutative by checking that all
the inner squares commute.  The lower right square is commutative as a
the result of applying $\Phi$ followed by $\mathcal{H}^{0}$ to the
descent diagram for $(\{ \mathcal{S}_{i} \}, \{\kappa_{ij} \})$.  This
descent diagram can be found as diagram \ref{eqn:reenterprettwist}
where $A$ must be replaced by $\Upsilon(A)$.  The square to the left of the
lower right square is commutative using \ref{eqn:anotherthing}.  The
other squares of the diagram are clearly commutative.  Therefore the
functor $\widetilde{\phi}$ is well defined.  We define
$\widetilde{\Phi}$ to be the right derived functor $R\widetilde{\phi}$
of the left exact functor $\widetilde{\phi}$.
\[\widetilde{\Phi} = R\widetilde{\phi} :
\sff{D}^{b}_{c}({}_{(M, T, \zeta)} Y_{g}, -1) 
\to \sff{D}^{b}_{c}({}_{(L, \theta, \eta)} X_{f},-1) \]

\subsection{Twisted Derived Equivalences} \label{ssec:twistderive}

In this section, we prove that the functor $\widetilde{\Phi}$ defined
in Lemma \ref{lem:newfunctor} is actually an equivalence.  We start by
establishing some relationships between the pushforward of sheaves from
an atlas and the relative Fourier-Mukai transform.  Consider the
diagram

\[
\xymatrix{
& \mathcal{U} \times_{U} \mathcal{V}  \ar[dl]_-{\tilde{\rho}_{U}}  
\ar[dr]^-{\tilde{\pi}_{U}}
\ar[dd]_-{r}& \\
\mathcal{U} \ar[dd]_-{a} 
& & \mathcal{V} \ar[dd]_-{b} \\
& X \times_{B} Y \ar[dl]_-{\tilde{\rho}} \ar[dr]^-{\tilde{\pi}}
 & \\
X \ar[dr]_-{\pi} & & Y
\ar[dl]^-{\rho} \\
& B &
}
\]
Then we have a natural equivalence of functors 
$\text{Mod}(\mathcal{V}) \to \text{Mod}(X)$ given by
\begin{equation}  \label{eqn:FMpushunderived}
\phi \circ b_{*} \cong a_{*} \circ \phi_{U}.
\end{equation}
This implies that we have a natural equivalence of functors 
$\sff{D}(\mathcal{V}) \to \sff{D}(X)$ 
\begin{equation}  \label{eqn:FMpush}
\Phi \circ Rb_{*} \cong Ra_{*} \circ \Phi_{U}
\end{equation}
Indeed note that for any object $\mathcal{S}$ in
$\sff{D}(\mathcal{V})$ we have 
\[\Phi(R b_{*} \mathcal{S}) =
R \tilde{\rho}_{*}(\tilde{\pi}^{*}
(Rb_{*} \mathcal{S}) \otimes \mathcal{P}) \cong 
R \tilde{\rho}_{*}((Rr_{*}\tilde{\pi}_{U}^{*}\mathcal{S})\otimes \mathcal{P})
\cong R \tilde{\rho}_{*}(Rr_{*}(\tilde{\pi}_{U}^{*}\mathcal{S} 
\otimes \mathcal{P}_{U}))
\]
and 
\[R \tilde{\rho}_{*}(Rr_{*}(\tilde{\pi}_{U}^{*}\mathcal{S} 
\otimes \mathcal{P}_{U})) 
\cong Ra_{*}(R\tilde{\rho}_{U*}(\tilde{\pi}_{U}^{*} \mathcal{S} 
\otimes \mathcal{P}_{U})) \cong Ra_{*} \Phi_{U}(\mathcal{S}).
\]

Recall that our choice of presentations gives us maps
$\mathfrak{a} = (a,L): \mathcal{U} \to \mathfrak{X}$, and
$\mathfrak{b} = (b,M): \mathcal{V} \to \mathfrak{Y}$, fitting into the
following commutative diagrams

\[
\xymatrix{ & & \mathfrak{Y} \ar[d]^-{c_{\mathfrak{Y}}} \\
Y & \ar[l]^-{b} \mathcal{V} \ar[r]^-{b_{g}} \ar[ur]^-{\mathfrak{b}} & Y_{g}
}
\]
\noindent
and

\[
\xymatrix{ & & \mathfrak{X} \ar[d]^-{c_{\mathfrak{X}}} \\
X & \ar[l]^-{a} \mathcal{U} \ar[r]^-{a_{f}} \ar[ur]^-{\mathfrak{a}} & X_{f}
}
\]
The pushforward functor $\mathfrak{a}_{*,-1}$ is defined as 
\[\mathfrak{a}_{*,-1}: \text{Mod}(\mathcal{U}) 
\to \text{Mod}({}_{(L, \theta, \eta)} X_{f}, -1)
\]
\[\mathfrak{a}_{*,-1}(\mathcal{S}) = 
(p_{1*}(A p_{0}^{*}\mathcal{S}), p_{12*}
(\btheta(p_{01}^{*}p_{0}^{*}\mathcal{S}))).
\]
We use the expression
$p_{12*}(\btheta(p_{01}^{*}p_{0}^{*}\mathcal{S}))$ here as shorthand.  
Indeed recall that we have a natural transformation $\btheta:
p_{12}^{*}A \circ p_{01}^{*}A \to p_{02}^{*}A$ and therefore an
isomorphism of sheaves 
\[\btheta(p_{01}^{*}p_{0}^{*}\mathcal{S}): 
(p_{12}^{*}A)(p_{01}^{*}(A(p_{0}^{*}\mathcal{S})) \to
(p_{02}^{*}A)(p_{01}^{*}p_{0}^{*}\mathcal{S}).
\]
We can use this to get a chain of isomorphisms 
\[Ap_{0}^{*}p_{1*}(Ap_{0}^{*}\mathcal{S}) 
\to p_{12*}((p_{12}^{*}A)(p_{01}^{*}(A(p_{0}^{*}\mathcal{S}))) 
\to p_{12*}((p_{02}^{*}A)(p_{01}^{*}p_{0}^{*}\mathcal{S})) 
\to p_{1}^{*}p_{1*}(A p_{0}^{*}\mathcal{S})
\]
as needed (see remark \ref{rem:reints}).
The associated derived functor is given by
\[R\mathfrak{a}_{*,-1}: \sff{D}(\mathcal{U}) 
\to  \sff{D}({}_{(L, \theta, \eta)} X_{f}, -1)
\]
\[R\mathfrak{a}_{*,-1}(\mathcal{S}) = 
(Rp_{1*}(A p_{0}^{*}\mathcal{S}), Rp_{12*}
(\btheta(p_{01}^{*}p_{0}^{*}\mathcal{S}))).
\]
Similarly, the pushforward functor $\mathfrak{b}_{*,-1}$ is defined as 
\[\mathfrak{b}_{*,-1}: \text{Mod}(\mathcal{V}) 
\to \text{Mod}({}_{(M, T, \zeta)} Y_{g}, -1)
\]
\[\mathfrak{b}_{*,-1}(\mathcal{T}) = 
(q_{1*}(\Upsilon (A) q_{0}^{*}\mathcal{T}), q_{12*}
(\bT(q_{01}^{*}q_{0}^{*}\mathcal{T}))).
\]
and its associated derived functor is given by 
\[R\mathfrak{b}_{*,-1}: \sff{D}(\mathcal{V}) 
\to \sff{D}({}_{(M, T, \zeta)} Y_{g}, -1)
\]
\[R\mathfrak{b}_{*,-1}(\mathcal{T}) = 
(Rq_{1*}(\Upsilon (A) q_{0}^{*}\mathcal{T}), Rq_{12*}
(\bT(q_{01}^{*}q_{0}^{*}\mathcal{T}))).
\]


Analogously to \ref{eqn:FMpushunderived}, we have a 
natural equivalence of functors 
$\text{Mod}(\mathcal{V}) \to \text{Mod}({}_{(L, \theta, \eta)} X_{f})$ 
given by 
\[\tilde{\phi} \circ \mathfrak{b}_{*, -1} 
\cong \mathfrak{a}_{*, -1} \circ \phi_{U}.
\]
which induces the natural equivalence of functors
$\sff{D}(\mathcal{V}) \to \sff{D}({}_{(L, \theta, \eta)} X_{f})$ given as 
\begin{equation} \label{eqn:twistFMpush}
\widetilde{\Phi} \circ R\mathfrak{b}_{*, -1} \cong 
R\mathfrak{a}_{*, -1} \circ \Phi_{U}.
\end{equation}
To see this, first note that applying the above to

\[
\xymatrix{
& \mathcal{U}^{1} \times_{U^{1}} \mathcal{V}^{1}  
\ar[dl]_-{\tilde{\rho}_{U^{1}}}  
\ar[dr]^-{\tilde{\pi}_{U^{1}}}
\ar[dd]_-{r_{U}}& \\
\mathcal{U}^{1} \ar[dd]_-{p_{1}} 
& & \mathcal{V}^{1} \ar[dd]_-{q_{1}} \\
& \mathcal{U}^{} \times_{U^{}} \mathcal{V}^{} 
\ar[dl]_-{\tilde{\rho}_{U}} \ar[dr]^-{\tilde{\pi}_{U}}
 & \\
\mathcal{U} \ar[dr]_-{\pi_{U}} & & \mathcal{V}
\ar[dl]^-{\rho_{U}} \\
& U &
}
\]
gives us
\[\phi_{U} \circ q_{1*} \cong p_{1*} \circ \phi_{U^{1}}
\]
as functors 
$\text{Mod}(\mathcal{V}^{1}) \to \text{Mod}(\mathcal{U})$.  Therefore, 
for any object $\mathcal{S}$ in $\text{Mod}(\mathcal{V})$, we have 
\[\widetilde{\phi}  (\mathfrak{b}_{*, -1} \mathcal{S}) = 
\widetilde{\phi} (q_{1*} (\Upsilon(A)(q_{0}^{*} \mathcal{S})), 
q_{12*}(T(q_{01}^{*}q_{0}^{*}\mathcal{S})))
\]
\[ \cong 
(\widetilde{\phi}(q_{1*} (\Upsilon(A)(q_{0}^{*} \mathcal{S})), 
\widetilde{\phi}(q_{12*}(T(q_{01}^{*}q_{0}^{*}\mathcal{S})))) 
\]
But recall that 
\[\widetilde{\phi}(q_{1*} (\Upsilon(A)(q_{0}^{*} \mathcal{S})) = 
\mathcal{H}^{0}(\Phi_{U} (q_{1*} (\Upsilon(A)(q_{0}^{*} \mathcal{S}))))
\cong \phi_{U}(q_{1*}(\Upsilon(A)q_{0}^{*} \mathcal{S}))
\]
however,
\[\phi_{U} (q_{1*} (\Upsilon(A)(q_{0}^{*} \mathcal{S}))) 
\cong p_{1*} (\phi_{U^{1}} (\Upsilon(A)(q_{0}^{*} \mathcal{S})))
\]
and 
\[
p_{1*} (\phi_{U^{1}} (\Upsilon(A)(q_{0}^{*} \mathcal{S}))) 
\cong p_{1*} (A (\phi_{U^{1}} q_{0}^{*} \mathcal{S})) 
\cong p_{1*} (A p_{0}^{*} (\phi_{U}  \mathcal{S}))  
\]
Taking these isomorphisms into account, we also have 
\[\widetilde{\phi}(q_{12*}(\bT(q_{01}^{*}q_{0}^{*}\mathcal{S}))) 
\cong p_{12*}(\btheta(p_{01}^{*}p_{0}^{*}(\phi_{U} \mathcal{S}))),
\]
this is essentially the content of the diagram \ref{eqn:bigupsilongerbe}.
Therefore, $\widetilde{\phi}(\mathfrak{b}_{*,-1} (\mathcal{S})) \cong
\mathfrak{a}_{*,-1}(\phi_{U} (\mathcal{S}))$, and taking derived
functors, we have demonstrated the isomorphism \ref{eqn:twistFMpush}.

Let $c: \mathfrak{S} \to S$ be a gerbe on a compact, connected,
complex manifold $S$.  We will need two facts about the triangulated
categories $D^{b}_{c}(\mathfrak{S},k)$, both of which were proven in
$\cite{DP}$.
\begin{itemize}
\item  $D^{b}_{c}(\mathfrak{S},k)$ admits 
a Serre functor defined by 
\[\bS_{\mathfrak{S},k}: D^{b}_{c}(\mathfrak{S},k) 
\to D^{b}_{c}(\mathfrak{S},k)
\]
\[a \mapsto a \otimes c^{*} \omega_{X}[n]
\]
\item
The structure sheaves of points on $S$ viewed as 
sheaves of weight $k$, define a spanning class for 
$D^{b}_{c}(\mathfrak{S},k)$.
\end{itemize}
We now recall a general criterion, due to Bondal-Orlov, and Bridgeland, for 
a functor $\bF : \mycal{A} \to \mycal{B}$ to be an 
equivalence of linear triangulated categories.

{\bf Theorem \cite{bondal-orlov-flops, bridgeland,bkr}} 
{\em Assume that $\mycal{A}$ and
$\mycal{B}$ have Serre functors $\bS_{\!\mycal{A}}$,
$\bS_{\!\mycal{B}}$, that $\mycal{A} \neq 0$, $\mycal{A}$ has a 
spanning class $C$, $\mycal{B}$ is
indecomposable, and that $\bF : \mycal{A} \to \mycal{B}$ has left
and right adjoints. Then $\bF$ is an equivalence if 
\[\bF : \op{Hom}_{\mycal{A}}^{i}(x_{1},x_{2}) \widetilde{\to}
\op{Hom}_{\mycal{B}}^{i}(\bF x_{1},\bF x_{2}), \quad \text{for
all } i \in {\mathbb Z}, x_{1}, x_{2} \in C
\]
and it intertwines the Serre functors: 
\[ \bF\circ \bS_{\!\mycal{A}}(x) =
\bS_{\!\mycal{B}}\circ \bF(x)
\]
on all elements $x \in
C$ in the spanning class. }

\begin{theo} \label{theo:main}
Let $X$ and $Y$ be compact, connected, complex manifolds and $\pi: X
\to B$ and $\rho:Y \to B$ be maps of complex analytic spaces, with
$\rho$ flat and $\mathcal{P} \in \text{Coh}(X \times_{B} Y)$ which is
flat over $Y$, implement a derived equivalence $\Phi: \sff{D}^{b}_{c}(Y) \to
\sff{D}^{b}_{c}(X)$.  Let $\mathfrak{X}$ be a $\Phi$-compatible gerbe on a
twisted version of $X$.  Let $\mathfrak{Y}$ be a $\Phi$-dual gerbe to
$\mathfrak{X}$ (so $\mathfrak{Y}$ lives over a twisted version of
$Y$).  Then we have an equivalence of categories
\[\sff{D}^{b}_{c}(\mathfrak{Y}, -1) \to \sff{D}^{b}_{c}(\mathfrak{X}, -1)
\]
\end{theo}

{\bf Proof.}  First of all, note that our functor $\widetilde{\Phi}$ is
a right-derived functor, and hence has a left adjoint.  Since
$\sff{D}^{b}_{c}(\mathfrak{Y}, -1)$ and $\sff{D}^{b}_{c}(\mathfrak{X},
-1)$, both have Serre functors, $\widetilde{\Phi}$ therefore has a
right adjoint as well.  It is obvious that
$\sff{D}^{b}_{c}(\mathfrak{Y}, -1)$ is not zero, and that
$\sff{D}^{b}_{c}(\mathfrak{X}, -1)$ is indecomposable, and it has
already been noted \cite{DP} that the points give a spanning class of
$\sff{D}^{b}_{c}(\mathfrak{Y}, -1)$.  

Therefore, we simply need to show that $\widetilde{\Phi}$ is
orthogonal on this spanning class, and that it intertwines the Serre
functors on this spanning class.

We claim that $\widetilde{\Phi}$ is orthogonal on the spanning class
consisting of points in $Y_{g}$.  Let $v_{1}$ and $v_{2}$ be any two
points in $\mathcal{V}$.  Then we know that
\[\text{Hom}^{i}_{D^{b}_{c}(Y)}
(Ra_{*}\mathbb{C}_{v_{1}},Ra_{*}\mathbb{C}_{v_{2}}) 
\to \text{Hom}^{i}_{D^{b}_{c}(X)}(\Phi Ra_{*}
\mathbb{C}_{v_{1}},\Phi Ra_{*}\mathbb{C}_{v_{2}}) \cong
\]
\[
\cong \text{Hom}^{i}_{D^{b}_{c}(X)}(Rb_{*} \Phi_{U}
\mathbb{C}_{v_{1}}, Rb_{*} \Phi_{U} \mathbb{C}_{v_{2}}) 
\]
is an isomorphism, where we have used the isomorphism \ref{eqn:FMpush}.
Now the commutative diagram 
\[
\xymatrix{ \sff{D}^{b}_{c}({}_{(M, T, \zeta)} Y_{g}, -1) 
\ar[r]^-{\widetilde{\Phi}} \ar[d]^-{} & 
\sff{D}^{b}_{c}({}_{(L, \theta, \eta)} X_{f}, -1) \ar[d]^-{} \\
\sff{D}^{b}(\mathcal{V})  \ar[r]^-{\Phi_{U}} 
& \sff{D}^{b}(\mathcal{U}) \\
\sff{D}^{b}_{c}(Y) \ar[r]^-{\Phi} \ar[u]^-{} & \sff{D}^{b}_{c}(X) \ar[u]^-{}
}
\]
induces the following commutative diagram where the vertical arrows
are isomorphisms, with inverses given by the various pushforward
functors, see \ref{eqn:FMpush} and \ref{eqn:twistFMpush}.
\[
\xymatrix{ \text{Hom}^{i}_{\sff{D}^{b}_{c}
(\mathfrak{Y}, -1)}(R\mathfrak{b}_{*,-1}\mathbb{C}_{v_{1}}, 
R\mathfrak{b}_{*,-1}\mathbb{C}_{v_{2}})
\ar[r]^-{\widetilde{\Phi}} \ar[d]^-{\mathfrak{b}_{-1}^{*}} & 
\text{Hom}^{i}_{\sff{D}^{b}_{c}(\mathfrak{X}, -1)}(\widetilde{\Phi} 
\mathfrak{b}_{*,-1} \mathbb{C}_{v_{1}}, 
\widetilde{\Phi} \mathfrak{b}_{*,-1} \mathbb{C}_{v_{2}}) 
\ar[d]^-{\mathfrak{a}_{-1}^{*}} \\
\text{Hom}^{i}_{\sff{D}^{b}(\mathcal{V})}
(\mathbb{C}_{v_{1}},\mathbb{C}_{v_{2}}) \ar[r]^-{\Phi_{U}} 
\ar[r]^-{}  & \text{Hom}^{i}_{\sff{D}^{b}(\mathcal{U})}
(\Phi_{U} \mathbb{C}_{v_{1}},\Phi_{U} \mathbb{C}_{v_{2}}) \\
\text{Hom}^{i}_{\sff{D}^{b}_{c}(Y)}(Rb_{*}\mathbb{C}_{v_{1}},
Rb_{*}\mathbb{C}_{v_{2}}) \ar[r]^-{\Phi} \ar[u]^-{b^{*}} & 
\text{Hom}^{i}_{\sff{D}^{b}_{c}(X)}(\Phi Rb_{*}\mathbb{C}_{v_{1}},
\Phi Rb_{*}\mathbb{C}_{v_{2}}) \ar[u]^-{a^{*}}
}
\]

Since $\Phi$ is an equivalence, the bottom is an isomorphism and so we
are done.  Also using the fact that $\Phi$ is an equivalence, we know
that $\Phi$ intertwines the Serre functors of $\sff{D}^{b}_{c}(Y)$ and
$\sff{D}^{b}_{c}(X)$.  In other words 
\[\Phi \circ \bS_{Y} \cong \bS_{X} \circ \Phi.
\]
Now pulling back by the cover $U \to B$ and using the property
\ref{eqn:locality} we see that for all sheaves $\mathcal{S}$ on
$\mathcal{V}$ we have 
\begin{equation}\label{eqn:pullback_intertwine}
\Phi_{U}((b^{*} \omega_{Y}) \otimes \mathcal{S} [n]) 
\cong (a^{*}\omega_{X}) \otimes \Phi_{U}(\mathcal{S})[n] 
\end{equation}
Now we claim that $\widetilde{\Phi}$ intertwines
the Serre functors of $\sff{D}^{b}_{c}(\mathfrak{Y},-1)$ and
$\sff{D}^{b}_{c}(\mathfrak{X},-1)$ on the spanning class.  Consider a
point $v$ of $\mathcal{V}$.  Observe there is a sheaf
$\mathcal{S} = \mathbb{C}_{v}$, with the properties $Rb_{*} \mathcal{S} \cong
\mathbb{C}_{b(v)}$, $Rb_{g*} \mathcal{S} \cong \mathbb{C}_{b_{g}(v)}$,
and $R\mathfrak{b}_{*, -1} \mathcal{S} \cong
\mathbb{C}_{b_{\mathfrak{g}}(v), -1}$.  Then if we let $\mathcal{T} =
\Phi_{U}(\mathcal{S})$, we have $Ra_{*} \mathcal{T} \cong \Phi
(\mathbb{C}_{b(v)})$, and $R\mathfrak{a}_{*,-1}\mathcal{T} \cong
\widetilde{\Phi} (\mathbb{C}_{b_{g}(v), -1})$.  This follows from
equations \ref{eqn:FMpush} and \ref{eqn:twistFMpush}.

Moving forward, we have,

\[\widetilde{\Phi}(\bS_{\mathfrak{Y},-1}(\mathbb{C}_{b_{g}(v), -1})) 
= \widetilde{\Phi}((c_{\mathfrak{Y}}^{*} \omega_{Y_{g}}) 
\otimes \mathbb{C}_{b_{g}(v), -1}[n])  \cong 
\widetilde{\Phi}(c_{\mathfrak{Y}}^{*} 
\omega_{Y_{g}}) \otimes R\mathfrak{b}_{*, -1}(\mathcal{S})[n]) 
\]
and
\[ \widetilde{\Phi}(c_{\mathfrak{Y}}^{*} 
\omega_{Y_{g}}) \otimes R\mathfrak{b}_{g*, -1}(\mathcal{S})[n]) \cong 
\widetilde{\Phi} R\mathfrak{b}_{*, -1}(b_{g}^{*} 
\omega_{Y_{g}} \otimes \mathcal{S}[n]) 
\]
and by \ref{eqn:twistFMpush} and the fact that $\Phi$ intertwines the
classical Serre functors (\ref{eqn:pullback_intertwine})
\[
\widetilde{\Phi} R\mathfrak{b}_{*, -1}(b_{g}^{*} \omega_{Y_{g}} 
\otimes \mathcal{S}[n]) 
\cong R\mathfrak{a}_{*, -1} \Phi_{U}((b^{*}\omega_{Y}) 
\otimes \mathcal{S} [n]) \cong 
R\mathfrak{a}_{*, -1} ((a^{*}\omega_{X}) 
\otimes \Phi_{U}( \mathcal{S}) [n]) 
\]
however
\[
R\mathfrak{a}_{*, -1} ((a^{*}\omega_{X}) 
\otimes \Phi_{U}( \mathcal{S}) [n]) \cong
R\mathfrak{a}_{f*, -1} ((a_{f}^{*}\omega_{X_{f}}) \otimes \mathcal{T} [n]) 
\cong (c_{\mathfrak{X}}^{*} \omega_{X_{f}}) 
\otimes R\mathfrak{a}_{*,-1} \mathcal{T} [n]
\]
and finally 
\[ (c_{\mathfrak{X}}^{*} \omega_{X_{f}}) 
\otimes R\mathfrak{a}_{*,-1} \mathcal{T} [n]
\cong \bS_{\mathfrak{X},-1} (\widetilde{\Phi} (\mathbb{C}_{b_{g}(v), -1}))
\]
Therefore, the Serre functors are intertwined on the spanning class by
$\widetilde{\Phi}$, so we have checked all the criteria and hence
$\widetilde{\Phi}$ is an equivalence.

\hfill $\Box$

\subsection{Alternative Method}

In this subsection, we give another method for proving the
equivalence.  This method applies in different situations, allowing us
to drop various coherence, smoothness, and compactness assumptions.
We were convinced of the possibility/usefulness of such a method
due to comments of D. Arinkin who kindly provided his insights
to us after a reading of the thesis.  In the above, we have given an
assignment $\mathbb{A}$ which takes a pair of shifted sheaves
$\mathcal{P}, \mathcal{Q} \in \sff{D}_{c}(X \times_{B} Y)$ giving inverse
equivalences, and a compatible gerbe presentation $\widetilde{X} =
{}_{(L,\theta, \eta)} X_{f}$ to an equivalence $\sff{D}(\widetilde{Y})
\to \sff{D}(\widetilde{X})$, where $\widetilde{Y}$ is the dual
presentation to $\widetilde{X}$ defined by the pair $(\mathcal{P},
\mathcal{Q})$.  We encode this assignment as
\[(\widetilde{X}; \mathcal{P}, \mathcal{Q}) \leadsto 
[\mathbb{A}(\widetilde{X}; \mathcal{P}, \mathcal{Q}): 
\sff{D}(\widetilde{Y}(\widetilde{X}; \mathcal{P}, \mathcal{Q}),-1) 
\to \sff{D}(\widetilde{X},-1)]
\]
Observe that this assignment has the following properties:

\begin{itemize}
\item Given isomorphisms $\mathcal{P} \cong
\mathcal{P}'$, and $\mathcal{Q} \cong \mathcal{Q}'$, the following
diagram is commutative
\[\xymatrixcolsep{5pc}
\xymatrix{\sff{D}(\widetilde{Y}(\widetilde{X}; \mathcal{P}, \mathcal{Q}),-1) 
\ar@<0ex>[r]^-{\mathbb{A}(\widetilde{X}; \mathcal{P}, \mathcal{Q})} 
\ar@<0ex>[d]^-{}
& \sff{D}(\widetilde{X},-1) \\
\sff{D}(\widetilde{Y}(\widetilde{X}; \mathcal{P}', \mathcal{Q}'),-1) 
\ar@<0ex>[ur]_-{\mathbb{A}(\widetilde{X}; \mathcal{P}', \mathcal{Q}')}&
}
\]
\item
Suppose we are given shifted coherent sheaves $\mathcal{P}_{1},
\mathcal{Q}_{1}$ on $X \times_{B}Y$, and $\mathcal{P}_{2},
\mathcal{Q}_{2}$ on $Y \times_{B}Z$ giving inverse equivalences, and a
$(\mathcal{P}_{1},\mathcal{Q}_{1})$ compatible presentation
$\widetilde{X}$.  Let $\widetilde{Y}$ denote
$\widetilde{Y}(\widetilde{X};\mathcal{P}_{1},\mathcal{Q}_{1})$, and
suppose that $\widetilde{Y}$ is $(\mathcal{P}_{2},\mathcal{Q}_{2})$
compatible.  Let $\widetilde{Z}$ denote
$\widetilde{Z}(\widetilde{Y};\mathcal{P}_{2},\mathcal{Q}_{2})$.  Then
if $\mathcal{R} = \mathcal{P}_{1} * \mathcal{P}_{2}$ and $\mathcal{S}
= \mathcal{Q}_{2} * \mathcal{Q}_{1}$ are shifted sheaves, and we let
\[\widetilde{Z}'= \widetilde{Z}'(\widetilde{X};\mathcal{R},\mathcal{S})
\cong \widetilde{Z},
\]
the inverse equivalences given by $\mathcal{R}$ and $\mathcal{S}$ fit
into the following commutative diagram:
\[\xymatrixcolsep{5pc}
\xymatrix{\sff{D}(\widetilde{Z},-1)
\ar@<0ex>[r]^-{\mathbb{A}(\widetilde{Y}; \mathcal{P}_{2},
\mathcal{Q}_{2})} \ar@<0ex>[d]^-{}
& \sff{D}(\widetilde{Y},-1)
\ar@<0ex>[d]^-{\mathbb{A}(\widetilde{X}; \mathcal{P}_{1},
\mathcal{Q}_{1})} \\
\sff{D}(\widetilde{Z}',-1)  
\ar@<0ex>[r]^-{\mathbb{A}(\widetilde{X}; \mathcal{R},
\mathcal{S})} & \sff{D}(\widetilde{X},-1) 
}
\]
\item
If $Y=X$, $\mathcal{P} = \mathcal{Q} = \mathcal{O}_{\Delta}$, then
$\widetilde{Y}(\widetilde{X}; \mathcal{P} , \mathcal{Q}) \cong
\widetilde{X}$ and the resulting functor
$\mathbb{A}(\widetilde{X};\mathcal{O}_{\Delta},
\mathcal{O}_{\Delta}):\sff{D}(\widetilde{X}) \to
\sff{D}(\widetilde{X})$ is naturally equivalent to the identity.
\end{itemize}
\noindent
These properties immediately imply the following:
\begin{theo} \label{theo:main2}
Let $X, Y, B$ be complex analytic spaces and $\pi: X \to B$ and
$\rho:Y \to B$ be compact, flat maps.  Let $\mathcal{P}, \mathcal{Q}$
be shifted coherent sheaves on $X \times_{B} Y$, flat over both
factors, which implement inverse equivalences $\Phi:\sff{D}(Y) \to
\sff{D}(X)$, and $\Psi:\sff{D}(X) \to \sff{D}(Y)$, and hence satisfy
$\mathcal{P}*\mathcal{Q} \cong \mathcal{O}_{\Delta}$ and
$\mathcal{Q}*\mathcal{P} \cong \mathcal{O}_{\Delta}$.  Let
$\mathfrak{X}$ be a gerbe on a twisted version of $X$, compatible with
these equivalences.  Let $\mathfrak{Y}$ be the dual gerbe to
$\mathfrak{X}$, with respect to $\mathcal{P}$ and $\mathcal{Q}$ (so
$\mathfrak{Y}$ lives over a twisted version of $Y$).  Then we have an
equivalences of categories
\[\sff{D}^{*}(\mathfrak{Y}, -1) \to \sff{D}^{*}(\mathfrak{X}, -1)
\]
for $*= \emptyset, b$
and
\[\sff{D}^{*}_{c}(\mathfrak{Y}, -1) \to \sff{D}^{*}_{c}(\mathfrak{X}, -1)
\]
for $*= b, -$.
\end{theo}

{\bf Proof.}  By the above properties of $\mathbb{A}$, the functors
associated to $\mathcal{P}$ and $\mathcal{Q}$ compose to a functors
equivalent to the identity, since $\mathcal{P}$ and $\mathcal{Q}$
convolve to sheaves isomorphic to the diagonals.

\hfill $\Box$
\begin{remark}\label{rem:quasicoherent}
Using the definition of quasi-coherence found in \cite{BBP}, we expect
that one would also get equivalences
\[\sff{D}^{*}_{qc}(\mathfrak{Y}, -1) \to \sff{D}^{*}_{qc}(\mathfrak{X}, -1)
\]
for $*= b, -$.
\end{remark}

\section{An Application}
In this chapter we provide an application.  We would like to emphasize
that the geometric limitations we impose in the following are {\it
only necessary} in order to allow for the setup of Donagi and Pantev to be
reproduced via our setup.  The situations that we will consider in the
following are not the most general applications that one could imagine
of our main theorem, {\it even in the context of elliptic fibrations}.

In the first section, we look at an analogue of Donagi and
Pantev's setup for complex torus fibrations, and prove a general
corollary of our main theorem.  In the second section, we apply this
corollary to prove their conjectures.
\subsection{Complex Torus Fibrations} \label{ssec:cmplx_tor_fib} 
Consider a complex torus fibration $X \to B$ with section $\sigma:B
\to X$.  By this we just mean that away from a closed co-dimension one
analytic subset of $B$, the fiber of $\pi$ is a complex torus.  Away
from the singular fibers this is an analytic group space over $B$,
using the section as an identity.  Therefore over the complement of
the singular fibers, the translation by a section (relative to
$\sigma$) defines an automorphism of the map $X \to B$.
\begin{defi}
Let us call a complex torus fibration with section $X \to B$ {\it
reasonable} if $X$ is a compact, connected complex manifold, the map
to $B$ is flat and if the translation by any local section $U \to X$
extends uniquely to an automorphism of the map $X \times_{B} U \to U$.
\end{defi}
Let $\pi: X \to B$ and $\rho:Y \to B$ be reasonable complex torus
fibrations with sections over a common base $B$.  We use the section
to give a group structure to the sheaves of sections.  A sheaf $S$ on
$X \times_{B} Y$ is called a {\it bi-extension} if for any local
sections $\sigma: U \to X$ and $\psi: U \to Y$ we have that
$S|_{\sigma \times \rho^{-1}(U)}$ and $S|_{\pi^{-1}(U) \times \psi}$
are line bundles and we have isomorphisms
\[(p_{1} + p_{2}, p_{3})^{*}S \cong p_{13}^{*} S \otimes p_{23}^{*} S
\]
on $X \times_{U} \sigma \times_{U} Y$ and isomorphisms
\[(p_{1}, p_{2}+ p_{3})^{*}S \cong p_{12}^{*} S \otimes p_{13}^{*} S
\]
on $X \times_{U} Y \times_{U} \psi$ and these isomorphisms satisfy
certain natural compatibilities.  In particular notice this implies
that for small open sets $U \subset B$ and sections $\sigma$ of $X$
over $U$, and $\psi$ of $Y$ over $U$ that if $f$ (and $g$) represent
the translation action of $\sigma$ on $X$ ($\psi$ on $Y$) we have
\[f^{*} (S|_{\pi^{-1}(U) \times \psi}) \cong S|_{\pi^{-1}(U) \times \psi} 
\otimes (\pi^{*} ((\sigma, \psi)^{*} S|_{\pi^{-1}(U) \times \rho^{-1}(U)})) 
\]
Now if $U$ has no non-trivial line bundles, then the line bundle on
$U$ given by $((\sigma, \psi)^{*} S|_{\pi^{-1}(U) \times
\rho^{-1}(U)})$ can be trivialized, so we get an isomorphism
\begin{equation}\label{eqn-trans_inv_on_X}
f^{*} (S|_{\pi^{-1}(U) \times \psi}) \cong S|_{\pi^{-1}(U) \times \psi}.
\end{equation}
Similarly there exists an isomorphism, 
\begin{equation}\label{eqn-trans_inv_on_Y}
g^{*} (S|_{\sigma \times \rho^{-1}(U)}) 
\cong S|_{\sigma \times \rho^{-1}(U)}.
\end{equation}
Furthermore the bi-extension structure gives us isomorphisms 
\begin{equation}\label{eqn-X_tran_Y_tens}
(f \times 1)^{*} (S|_{\pi^{-1}(U) \times \rho^{-1}(U)}) 
\cong S \otimes \tilde{\pi}^{*} (S|_{\sigma \times \rho^{-1}(U)})
\end{equation}
and 
\begin{equation}\label{eqn-Y_tran_X_tens}
(1 \times g)^{*} (S|_{\pi^{-1}(U) \times \rho^{-1}(U)}) 
\cong S \otimes \tilde{\rho}^{*} (S|_{\pi^{-1}(U) \times \psi}).
\end{equation}

Suppose that $\pi: X \to B$ and $\rho: Y \to B$ are reasonable complex
torus fibrations with section over the same base $B$ and let us denote
by $\mathcal{X}$ and $\mathcal{Y}$ the sheaves of sections of $X$ and
$Y$ respectively.  By the above injections of sheaves, the sheaves
$\mathcal{X}$ and $\mathcal{Y}$ become sub-sheaves of groups in
$Aut_{X/B}$ and $Aut_{Y/B}$ respectively.  We
assume that there is a coherent Poincar\'{e} sheaf $\mathcal{P}$
living on $X \times_{B} Y$, flat over $Y$ that implements an
equivalence of derived categories $\Phi:\sff{D}^{b}_{c}(Y) \to
\sff{D}^{b}_{c}(X)$, and is a bi-extension.  Consider two elements
$\alpha \in H^{1}(B, \mathcal{Y})$ and $\beta \in H^{1}(B,
\mathcal{X})$.  Let $f \in \check{Z}^{1}(B, Aut_{X/B})$ and $g \in
\check{Z}^{1}(B, Aut_{Y/B})$ be the automorphisms given by translation
by some representatives for $\beta$ and $\alpha$ respectively.  We now
define a map
\[\mathcal{X} \to R^{1} \rho_{\alpha *} \mathcal{O}^{\times}
\]
Indeed we can send a local section $\gamma_{i}$ over $U_{i}$ to the
equivalence class 
\[[\mathcal{P}|_{\gamma_{i} \times \rho^{-1}(U_{i})}]
\in H^{1}(\rho^{-1}(U_{i}),\mathcal{O}^{\times}).
\]
Let $g_{ij}$ be the automorphism of $\rho^{-1}(U_{ij})$ corresponding
to the section $\alpha_{ij}$.  Now, because of equation
\ref{eqn-trans_inv_on_Y} (following from the reasonable nature of the
torus fibrations and the bi-extension property of $\mathcal{P}$) we
have, for sections $\gamma_{i}$ and $\gamma_{j}$ agreeing on the
overlap
\[g_{ij}^{*} [\mathcal{P}|_{\gamma_{i} \times \rho^{-1}(U_{i})}] = 
[\mathcal{P}|_{\gamma_{j} \times \rho^{-1}(U_{i})}].
\]
in $H^{1}(\rho^{-1}(U_{ij}), \mathcal{O}^{\times})$.  This shows the
map that we have described is well defined.  Similarly, we have a map
\[\mathcal{Y} \to R^{1} \pi_{\beta *} \mathcal{O}^{\times}.
\]
By taking the induced map on cohomology, we denote by
$S_{\alpha}(\beta)$ the image of $\beta$ under the map
\[
H^{1}(B, \mathcal{X}) \to H^{1}(B, R^{1} \rho_{\alpha *}
\mathcal{O}^{\times}).
\]
Similarly, we denote by $S_{\beta}(\alpha)$ the image of $\alpha$
under the map
\[
H^{1}(B, \mathcal{Y}) \to H^{1}(B, R^{1} \pi_{\beta *}
\mathcal{O}^{\times}).
\]  
Let us further assume that the images of $S_{\alpha}(\beta)$ and
$S_{\beta}(\alpha)$ (via the differentials $d_{g;2}$ and $d_{f;2}$
described in section \label{sec:classification}) vanish in $H^{3}(B,
\mathcal{O}^{\times})$ and that the image of $H^{2}(B,
\mathcal{O}^{\times})$ vanishes in both $H^{2}(X_{\beta},
\mathcal{O}^{\times})$ and $H^{2}(X_{\alpha}, \mathcal{O}^{\times})$.
Then, according to the spectral sequence for the fibrations we find
unique gerbes ${}_{\alpha}X_{\beta} \to X_{\beta}$ and
${}_{\beta}Y_{\alpha} \to Y_{\alpha}$ which correspond to
$S_{\beta}(\alpha)$ and $S_{\alpha}(\beta)$ respectively.  We then
have the following corollary of the main theorem.  For certain {\it
elliptic} fibrations $X=Y$ with some possible further restrictions on
the pair $(\beta, \alpha)$, this fact has already been proven by
Donagi and Pantev $\cite{DP}$.  For some other cases of $\beta$ and
$\alpha$ this was conjectured.  The following corollary proves the
conjectures (see e.g. conjecture 2.19) from $\cite{DP}$ and as well it
gives an alternative proof of the main results in the paper
$\cite{DP}$.
\begin{cor} \label{cor-DP}
In the above circumstance, we have 
\[\sff{D}^{b}_{c}({}_{-\beta}Y_{\alpha}, 1) \cong 
\sff{D}^{b}_{c}({}_{\alpha}X_{\beta}, 1)
\]
\end{cor}
{\bf Proof.}  

In order to prove this we simply need to find gerbe presentations
${}_{(L,\theta, \eta)} X_{f}$ of ${}_{-\alpha}X_{\beta}$ and ${}_{(M,T,
\zeta)} Y_{g}$ of ${}_{\beta}Y_{\alpha}$ which are dual in the sense
of the main theorem.  For then we will have 
\[\sff{D}^{b}_{c}({}_{-\beta}Y_{\alpha}, 1) \cong 
\sff{D}^{b}_{c}({}_{(M,T, \zeta)} Y_{g}, -1) \cong
\sff{D}^{b}_{c}({}_{(L,\theta, \eta)} X_{f}, -1) \cong
\sff{D}^{b}_{c}({}_{- \alpha}X_{\beta}, -1) \cong
\sff{D}^{b}_{c}({}_{\alpha}X_{\beta}, 1)  
\]
We do this simply by picking a reasonable choice
for ${}_{(L,\theta, \eta)} X_{f}$ and showing that its dual
presentation represents ${}_{\beta}Y_{\alpha}$.

We know that the gerbe ${}_{-\alpha}X_{\beta}$ comes from the
$E_{\infty}$ term of the Leray-Serre spectral sequence, given by
\[ [[\mathcal{P}|_{(\mathcal{U}^{1})_{f} \times - \alpha}]] \in 
\frac{\text{ker}(H^{1}((\mathcal{U}^{1})_{f}, \mathcal{O}^{\times}) 
\to H^{1}((\mathcal{U}^{2})_{f}, \mathcal{O}^{\times}))}
{\text{im}(H^{1}((\mathcal{U}^{0})_{f}, \mathcal{O}^{\times}) 
\to H^{1}((\mathcal{U}^{1})_{f}, \mathcal{O}^{\times}))} 
\cong H^{1}(B, R^{1}\pi_{\beta*} \mathcal{O}^{\times})
\]
In other words, $[[\mathcal{P}|_{(\mathcal{U}^{1})_{f} \times
- \alpha}]]$ goes to zero in $H^{3}(B, \mathcal{O}^{\times})$.  The
double brackets refer to the fact that
$[[\mathcal{P}|_{(\mathcal{U}^{1})_{f} \times - \alpha}]]$ comes from the
element $[\mathcal{P}|_{(\mathcal{U}^{1})_{f} \times - \alpha}] \in
H^{1}((\mathcal{U}^{1})_{f}, \mathcal{O}^{\times})$.  Concretely, we
can see that $[\mathcal{P}|_{(\mathcal{U}^{1})_{f} \times - \alpha}] \in
\text{ker}(H^{1}((\mathcal{U}^{1})_{f}, \mathcal{O}^{\times}) \to
H^{1}((\mathcal{U}^{2})_{f}, \mathcal{O}^{\times}))$ by the equation 
\[[\mathcal{P}|_{(\mathcal{U}^{1})_{f} \times 
- \alpha_{jk}}] \otimes f_{jk}^{*}[\mathcal{P}|_{(\mathcal{U}^{1})_{f} 
\times - \alpha_{ij}}] = [\mathcal{P}|_{(\mathcal{U}^{1})_{f} 
\times - \alpha_{jk}}] \otimes [\mathcal{P}|_{(\mathcal{U}^{1})_{f} 
\times - \alpha_{ij}}]
\]
\[= [\mathcal{P}|_{(\mathcal{U}^{1})_{f} 
\times (- \alpha_{ij} - \alpha_{jk})}]=[\mathcal{P}|_{(\mathcal{U}^{1})_{f} 
\times - \alpha_{ik}}] 
\]
where we have used equation \ref{eqn-trans_inv_on_X}.
Consider the natural filtration
\[0 \subset F^{2} 
\mathbb{H}^{2}(\check{C}^{\bullet} ((\mathcal{U}^{\bullet})_{f},
\mathcal{O}^{\times})) \subset F^{1}
\mathbb{H}^{2}(\check{C}^{\bullet} ((\mathcal{U}^{\bullet})_{f},
\mathcal{O}^{\times})) \subset \mathbb{H}^{2}(\check{C}^{\bullet}
((\mathcal{U}^{\bullet})_{f}, \mathcal{O}^{\times}))
\]
Here $F^{1} \mathbb{H}^{2}(\check{C}^{\bullet}
((\mathcal{U}^{\bullet})_{f}, \mathcal{O}^{\times}))$ are the classes
representable by an element whose component in
$\check{C}^{2}((\mathcal{U}^{0})_{f}, \mathcal{O}^{\times})$ is
trivial, and $F^{2} \mathbb{H}^{2}(\check{C}^{\bullet}
((\mathcal{U}^{\bullet})_{f}, \mathcal{O}^{\times}))$ are the classes
representable by an element whose components in
$\check{C}^{2}((\mathcal{U}^{0})_{f}, \mathcal{O}^{\times})$ and
$\check{C}^{1}((\mathcal{U}^{1})_{f}, \mathcal{O}^{\times})$ are
trivial.  The map 
\[F^{1} \mathbb{H}^{2}(\check{C}^{\bullet}
((\mathcal{U}^{\bullet})_{f}, \mathcal{O}^{\times})) \to E_{\infty}^{1,1}
\]
sends an element $(L,\theta) \in \check{C}^{1}((\mathcal{U}^{1})_{f},
\mathcal{O}^{\times}) \times \check{C}^{0}((\mathcal{U}^{2})_{f},
\mathcal{O}^{\times})$ to the image in $E_{\infty}^{1,1}$ of the
obvious modification $L' \in \check{C}^{1}((\mathcal{U}^{1})_{f},
\mathcal{O}^{\times})$ of $L$ given by choosing an element in
$\check{C}^{0}((\mathcal{U}^{1})_{f}, \mathcal{O}^{\times})$ which
trivializes $\theta$.  Notice, that we have just replaced $L$ by an
isomorphic line bundle $L'$, in other words, $[L] = [L'] \in
H^{1}((\mathcal{U}^{1})_{f}, \mathcal{O}^{\times})$.  In the other
direction consider the following set theoretic splitting of the surjection
\[F^{1} \mathbb{H}^{2}(\check{C}^{\bullet}
((\mathcal{U}^{\bullet})_{f}, \mathcal{O}^{\times})) \to
E^{1,1}_{\infty}.\] Given an element $[[L]] \in E^{1,1}_{\infty}$, we
know that $[L]$ goes to zero in $H^{1}((\mathcal{U}^{2})_{f},
\mathcal{O}^{\times})$ and therefore the image of $L$ in
$\check{C}^{1}((\mathcal{U}^{2})_{f}, \mathcal{O}^{\times})$ is
bounded by an element $\theta$ in
$\check{C}^{0}((\mathcal{U}^{2})_{f}, \mathcal{O}^{\times})$.  Since
$[[L]]$ survives to $E_{3}^{1,1} = E_{\infty}^{1,1}$, $\theta$ goes to
the trivial element in $\check{C}^{0}((\mathcal{U}^{3})_{f},
\mathcal{O}^{\times})$.  Therefore, the pair $(L, \theta)$ gives an
element in $F^{1} \mathbb{H}^{2}(\check{C}^{\bullet}
((\mathcal{U}^{\bullet})_{f}, \mathcal{O}^{\times}))$.  Because we are
assuming that $F^{2} \mathbb{H}^{2}(\check{C}^{\bullet}
((\mathcal{U}^{\bullet})_{f}, \mathcal{O}^{\times}))$ is trivial (the
pullback of a gerbes from the base is trivial), the inverse map (the
set theoretic splitting) is actually a group homomorphism and we have
written explicitly the two maps in the isomorphism $F^{1}
\mathbb{H}^{2}(\check{C}^{\bullet} ((\mathcal{U}^{\bullet})_{f},
\mathcal{O}^{\times})) \cong E_{\infty}^{1,1}$.

We now apply the above in the case that $L =
\mathcal{P}|_{(\mathcal{U}^{1})_{f} \times - \alpha}$, in other words
$L_{ij} =\mathcal{P}|_{(\mathcal{U}^{1})_{f} \times - \alpha_{ij}}$.  We
have canonical trivializations $\eta_{i}: \mathcal{O} \to L_{ii}$, and
therefore we have produced a presentation of a gerbe ${}_{(L, \theta,
\eta)} X_{f}$ on $X$ Now the gerbe presentation ${}_{(L,\theta, \eta)}
X_{f}$ corresponds, according to the proof of Lemma
\ref{lem:reenterpret}, to the element
\[ 
T_{\mathcal{P}|_{(\mathcal{U}^{1})_{f} 
\times -\alpha}}   \circ f^{*} \in \text{Aut}(\text{Coh}(\mathcal{U}^{1})).
\]
Pre-composing this by the Fourier-Mukai transform 
\[\Phi_{\mathcal{P}}: \sff{D}(\mathcal{V}^{1}) \to \sff{D}(\mathcal{U}^{1}) \]
yields
\[T_{\mathcal{P}|_{(\mathcal{U}^{1})_{f} 
\times -\alpha}}   \circ f^{*} \circ \Phi_{\mathcal{P}}
\]
Now using the equations \ref{eqn:Huy2} and \ref{eqn-X_tran_Y_tens} we have 
\[f^{*} \circ \Phi_{\mathcal{P}} \cong \Phi_{(f \times 1)^{*}\mathcal{P}}
\cong \Phi_{\mathcal{P} \otimes {\tilde{\pi}}^{*}(\mathcal{P}|_{\beta \times
\mathcal{V}^{1}})} \cong \Phi_{\mathcal{P}} \circ
T_{\mathcal{P}|_{\beta \times \mathcal{V}^{1}}}.
\]
Now we have, using equations \ref{eqn:Huy3} and
\ref{eqn-Y_tran_X_tens} we have
\[T_{\mathcal{P}|_{(\mathcal{U}^{1})_{f} 
\times - \alpha}} \circ \Phi_{\mathcal{P}} \cong
\Phi_{{\tilde{\rho}}^{*}(\mathcal{P}|_{(\mathcal{U}^{1})_{f} \times
- \alpha}) \otimes \mathcal{P}} \cong \Phi_{(1 \times g^{-1})^{*}
\mathcal{P}} \cong \Phi_{\mathcal{P}} \circ g^{*}.
\]
By putting all this together, we arrive at 
\[T_{\mathcal{P}|_{(\mathcal{U}^{1})_{f} 
\times - \alpha}} \circ f^{*} \circ \Phi_{\mathcal{P}} \cong
T_{\mathcal{P}|_{(\mathcal{U}^{1})_{f} \times - \alpha}} \circ
\Phi_{\mathcal{P}} \circ T_{\mathcal{P}|_{\beta \times
\mathcal{V}^{1}}} \cong \Phi_{\mathcal{P}} \circ g^{*} \circ
T_{\mathcal{P}|_{\beta \times \mathcal{V}^{1}}}.\] 
Therefore
\begin{equation} \label{eqn-compat}
T_{\mathcal{P}|_{(\mathcal{U}^{1})_{f} 
\times - \alpha}} \circ f^{*} \circ \Phi_{\mathcal{P}} \cong
\Phi_{\mathcal{P}} \circ T_{g^{*}(\mathcal{P}|_{\beta \times
\mathcal{V}^{1}})} \circ g^{*}.
\end{equation}
Equation \ref{eqn-compat} shows that the presentation ${}_{(L, \theta,
\eta)} X_{f}$ is $\Phi_{\mathcal{P}}$-compatible (see definition
\ref{defi:compatable}), and therefore the gerbe ${}_{-\alpha}
X_{\beta}$ is $\Phi_{\mathcal{P}}$-compatible and so we can calculate
the dual gerbe on a twisted version of $Y$.  Therefore, by the lemma
we have a dual gerbe, on the presentation $Y_{g}$ of $Y_{\alpha}$,
given by $M = g^{*}(\mathcal{P}|_{\beta \times \mathcal{V}^{1}}) \to
\mathcal{V}^{1}$, along with some isomorphisms $T_{ijk}:M_{jk} \otimes
g_{jk}^{*} M_{ij} \to M_{ik}$, and $\zeta_{i}: \mathcal{O} \to M_{ii}$
satisfying the diagrams of a gerbe presentation.  According to the
above discussion, the equivalence class of this gerbe on $Y_{\alpha}$
comes from the term in $E^{1,1}_{\infty}$ given by $[[M]] =
[[g^{*}(\mathcal{P}|_{\beta \times \mathcal{V}^{1}})]] =
[[\mathcal{P}|_{\beta \times \mathcal{V}^{1}}]]$.  This is the gerbe
on $Y_{\alpha}$ which comes from $\beta \in H^{1}(B, \mathcal{X})$ and
so we are done.

\hfill $\Box$

\begin{remark}\label{rem:more}
Clearly, we can also prove a more general statement by invoking
Theorem \ref{theo:main2}, we leave the recording of this statement to the
interested reader.
\end{remark}
\subsection{The Conjecture of Donagi and Pantev} 
\label{ssec:conjecture} 

In this section we use a special case of the corollary
\ref{cor-DP} and use it to reprove the main results in
$\cite{DP}$ as well as the following conjecture $\cite{DP}$ of Donagi
and Pantev, which was proven in $\cite{DP}$ in many important special
cases.
\begin{con} \label{con:DP}
Let $X$ be a complex manifold elliptically fibered with at worst $I_{1}$
fibers over a normal analytic variety $B$ such that $H^{2}(B,
\mathcal{O}^{\times}) = \{1 \}$.  Let $\alpha, \beta \in \TSh_{an}(X)$
be {\it complementary} elements (see \cite{DP} for the definition).
Then there exists an equivalence
\begin{equation} \label{eqn:whatwewant}
\sff{D}^{b}_{c}({}_{\beta} X_{\alpha}, -1) \cong 
\sff{D}^{b}_{c}({}_{\alpha} X_{\beta}, 1) 
\end{equation}
of the bounded derived categories of sheaves of weights $\pm 1$ on
${}_{\alpha} X_{\beta}$ and ${}_{\beta} X_{\alpha}$ respectively.
\end{con}
{\bf Proof}

Let $Y = X$ and let $\mathcal{X}$ be the sheaf of sections of $X$,
considered as a sheaf of groups via a choice of a section $\sigma$.
We have two flat morphisms, $\pi:X \to B$, and $\rho: X \to B$ giving
$X$ two different structures of a reasonable torus fibration over $B$.
We then get an isomorphism
\[\TSh_{an}(X) \cong H^{1}(B, \mathcal{X}).
\]
Consider the rank one divisorial sheaf on $X \times_{B} X$ defined as
\[\mathcal{P} = \mathcal{O}_{X \times_{B} X}(\Delta - \sigma \times_{B} X - 
X \times_{B} \sigma - \varpi^{*}N_{\sigma/X})
\]
where $\varpi: X \times_{B} X \to B$ is the natural projection.  It
gives an equivalence of derived categories $\Phi: \sff{D}^{b}_{c}(X)
\to \sff{D}^{b}_{c}(X)$ (see for example the paper of Bridgeland and
Maciocia \cite{BM}).  Furthermore, as explained in \cite{DP},
$\mathcal{P}$ is a bi-extension.  The fact that $\alpha$ and $\beta$
are complimentary (the triviality of the pairing
$<S_{\alpha}(\beta),S_{\beta}(\alpha)>$ from \cite{DP}) is equivalent
to the vanishing of the images of $S_{\alpha}(\beta)$ and
$S_{\beta}(\alpha)$ in $H^{3}(B, \mathcal{O}^{\times})$ as explained
in \cite{DP}.  Now an application of Corollary \ref{cor-DP}
immediately produces an isomorphism as in \ref{eqn:whatwewant} and
hence proves the conjecture \ref{con:DP} of Donagi and Pantev
\cite{DP}.

 \hfill $\Box$

\begin{remark}\label{rem:nodescusps}
In fact, since genus one fibrations with only nodes and cusps as
singularities are reasonable complex torus fibrations, we have
expanded the context in which their conjecture holds to this case.  Of
course, an even more general result for elliptic fibrations is given if
we use Theorem \ref{theo:main2}.
\end{remark}

\section{Conclusions and Speculations}

In this section, we comment on some possible future generalizations of
our main theorem.  After this, we indicate a few theorems that we
expect as consequences or analogues of the results of this thesis.
First of all, we expect to be able to dispense with the compactness
and smoothness assumptions of our varieties, as well as to be able to
extend our main result to the case where the fibrations are not flat,
and perhaps the Poincar\'{e} sheaf is a more general object in the
derived category of the fiber product.  Certainly, we could replace
our Serre functor with a suitable object in a singular situation.
Compactness assumptions can be dispensed with by considering derived
categories of compactly supported sheaves.  The key ingredient needed
for generalizations seems to be to find a diagram of isomorphisms of
sheaves or complexes of sheaves which induce the diagram of natural
equivalences \ref{eqn:anotherthing}.  In many cases this should be
possible using the constraints on the Poincar\'{e} object coming from
the fact that it implements an isomorphism.  Another strategy towards
finding an object on the fiber product of the two gerbes is to use the
technique of cohomological descent for gluing objects in the
differential graded categories: locally the object is just the
Poincar\'{e} object, and the diagram of natural equivalences
\ref{eqn:anotherthing} becomes a descent data diagram in the
differential graded category of sheaves on the fiber product.

Ideally, there would also be a presentation free description of the
twisting of derived equivalences.  A scheme to do precisely this was
suggested to us by D. Arinkin during a discussion of the original
version of the thesis.  It is an extension of his duality for group
stacks perspective, as found in the appendix to \cite{DP}.  We hope that
future results generalizing those of our main theorem will find
presentation free descriptions.

In the search for applications of our main theorem to fibrations where
the generic fiber has no group structure, a natural place to start is
$K3$ fibrations.  However, as we now explain, the most naive attempts
seem to lead back to scenarios which are already understood.  First,
we start with two $K3$ surfaces $W$ and $Z$ such that $Z$ is a
fine moduli space for stable sheaves of Mukai vector $(r, H, s)$ on
$W$.  If $L$ is any line bundle on $W$, then the map
\[T_{L}: \sff{D}^{b}_{c}(W) \to \sff{D}^{b}_{c}(W) 
\]
acts on the Mukai vector by $T_{L}(r, H, s) = (r, rc_{1}(L) +H,
s+(c_{1}(L), H) + \frac{1}{2}r(c_{1}(L)^{2})$.  Since $c_{1}(L)$ is
non-torsion, $r$ must then be zero in order that $T_{L}$ preserves the
moduli problem, and can be used to build gerbes with the derived
equivalence between $Z$ and $W$.  Unfortunately, this together with
the fact that $(r, H, s)^{2}=0$ implies that $H^{2} = 0$.  Now we can
find a non-trivial line bundle $M$ such that $c_{1}(M) = H$.  Then $M$
defines an elliptic fibration $W \to \mathbb{P}^{1}$.  The sheaves on
$W$ corresponding to points on $Z$ will all be supported on the fibers
of this elliptic fibration.  This argument seems to go through in the
relative case (where $(r,H,s)$ now becomes a relative Mukai vector)
and this shows that a $K3$ fibration $W \to B$ admitting gerbes
compatible with a derived equivalence (which are not pull-backs from
$B$) becomes an elliptic fibration over a $\mathbb{P}^{1}$ fibration
over $B$, and the derived equivalence respects this elliptic fibration
structure.  If it works, this reasoning will only be bi-rational, and
refining it to an actual geometric statement may lead to new
applications to $K3$ fibrations.  Also, note that if we replace our
$K3$ fibration with a, say, abelian surface fibration, then we could
have $r$ be non-zero, provided that $L$ was $r$-torsion, and hence we
expect a straightforward application of our main theorem in that case.

A perhaps more obvious application concerns the moduli spaces of
vector bundles on curves.  For concreteness, let $C$ be a genus $2$
hyperelliptic curve.  Then from the work of Desale and Ramanan
\cite{DeRa} we recall that the moduli space of (isomorphism classes of
stable) vector bundles of rank $2$ and fixed odd determinant of odd
degree is isomorphic to the intersection $X=Q_{1} \cap Q_{2}$ of two
quadric four-folds in $\mathbb{P}^{5}$.  Bondal and Orlov, in
\cite{bondal-orlov-flops} give an explanation in terms of derived
categories.  They produce a vector bundle $V$ on $C \times X$ that
implements a full and faithful embedding $\sff{D}^{b}_{c}(C) \to
\sff{D}^{b}_{c}(X)$.  Finally they prove that this morphism induces a
semiorthogonal decomposition:
\[\sff{D}^{b}_{c}(X) = 
\langle \mathcal{O}_{X}(-1), \mathcal{O}_{X}, \sff{D}^{b}_{c}(C) \rangle
\]
Now let $\pi:C \to B$, factorizing as $C \to P \to B$ be a smooth
compact family of non-singular hyperelliptic curves over a projective
bundle $P$ given by $6$ sections of $P$ over $B$.  Let $\mathfrak{C}$
be a gerbe on $C$, presentable by a collection of line bundles
$M_{ij}$ on $\pi^{-1}(U_{ij})$ which are $2-$torsion on each
hyperelliptic fiber.  Conjugating these line bundles by the
Fourier-Mukai transform (using a left adjoint) gives locally
automorphisms $f_{ij}$ of a natural family of quadrics $X$ inside a
$\mathbb{P}^{5}$ bundle over $B$.  By gluing we obtain a new family
of quadrics $\widetilde{X} \to B$.  A simple application of the main
theorem should give a full and faithful embedding
$\sff{D}^{b}_{c}(\mathfrak{C}, -1) \to
\sff{D}^{b}_{c}(\widetilde{X})$.  The line bundle
$\mathcal{O}_{X}(-1)$ will induce a line bundle
$\mathcal{O}_{\widetilde{X}}(-1)$ and the embedding will induce a
semiorthogonal decomposition:
\[\sff{D}^{b}_{c}(\widetilde{X}) = \langle \mathcal{O}_{\widetilde{X}}(-1), 
\mathcal{O}_{\widetilde{X}}, \sff{D}^{b}_{c}(\mathfrak{C},-1) \rangle.
\]
\subsection{Relation to Mirror Symmetry}
Finally, let us explain a simple example of the type of idea in
homological mirror symmetry that our main theorem might lead one to
conjecture.  We do not give a complete set of references here, but for
what we are going to discuss \cite{AZ} and the references therein are
relevant.  We would like to propose an extension of the usual ways
that gerbes are thought to enter the homological mirror symmetry story
\cite{Hitchin}, \cite{DP}.  However, our analysis applies to a
different type of gerbe.  Let $B_{1}=F_{1}=B_{2}=F_{2} = S^{1}$.  Let
$E_{1} = B_{1} \times F_{1}$, and $E_{2} = B_{2} \times F_{2}$ be
elliptic curves.  We assume that dualizing $F_{i}$ gives homological
mirror symmetry identifications $\sff{D}^{b}_{c}(E_{i}) \cong \sff{D}
Fuk(S_{i})$, which induce isomorphisms if we base change everything
via the inclusion of open sets in the base $B_{i}$.  Here $S_{i}=B_{i}
\times F_{i}^{\vee}$ are symplectic tori of real dimension $2$.  Let
$X= E_{1} \times E_{2}$ and $Y = S_{1} \times S_{2}$.  We also assume
we can take the "product" of equivalences and so we also have the
equivalence $\Phi: \sff{D} Fuk(Y) \to \sff{D}^{b}_{c}(X)$, which
induces isomorphisms if we base change everything via the inclusion of
open sets in $B_{1}$ (we will not care about the analogous property
for $B_{2}$).  Here $\sff{D} Fuk(Y)$ should include the coisotropic
branes of Kapustin and Orlov, see \cite{KO}.  Consider the projections
\[ \pi:X \to B_{1}, \ \ \ \
\rho:Y \to B_{1} \ \ \ \ \text{and} \ \ \ \ \mu: X \to E_{2}.
\]
Embed $E_{2}$ into projective space, this gives us a line bundle
$\mathcal{O}(1)$ on $E_{2}$.  Fix a class $\{n_{ij}\} \in
\check{Z}^{1}(B_{1}, \mathbb{Z})$ for some small open cover $\{U_{i}
\}$ of $B_{1}$.  Consider an $\mathcal{O}^{\times}-$gerbe
$\mathfrak{X}$ on $X$ presented by the collection of line bundles
$\{\mu^{*}\mathcal{O}(n_{ij}) \to \pi^{-1}(U_{ij}) \}$.  Now
conjugating the tensorization by $\mu^{*}\mathcal{O}(n_{ij})$ by the
(appropriately restricted) transformations $\Phi_{U_{ij}}$ yields
automorphisms $g^{n_{ij}}$ of $\rho^{-1}(U_{ij}) = U_{ij} \times
F_{1}^{\vee} \times S_{2}$.  These are symplectomorphisms of
$\rho^{-1}(U_{ij})$ which are constant on the $U_{ij} \times
F_{1}^{\vee}$ factor, and given by $n_{ij}$ compositions of a Dehn
twist $g$ on the factor $S_{2}$.  Using them, we can produce a new
symplectic $4$-manifold $\widetilde{Y}$.

Then, by analogy with the main theorem, we can conjecture that 
\[\sff{D}^{b}_{c}(\mathfrak{X}, -1) \cong \sff{D} Fuk(\widetilde{Y}).\]

Notice that the automorphisms destroy the second $S^{1}$ fibration
structure, but does nothing to the first one.  By contrast, if the
$\mathcal{O}(1)$ was replaced by a line bundle flat along the $F_{2}$
direction, we would expect to keep both $S^{1}$ fibration structures
on the symplectic side, fixing the first one and turning the second
into a non-trivial principal $S^{1}$ bundle.

Finally, we would like to mention a recent paper of Donagi and Pantev
on the duality between the Hitchin integrable systems associated to a
simple complex Lie group and its Langlands dual.  There, they prove a
derived equivalence between the bounded, compactly supported, coherent
derived categories of the commutative group stacks
${}^{L}\mathcal{H}iggs$ and $\mathcal{H}iggs$, over a space $B -
\Delta$.  Here $\Delta$ is a discriminant locus over which the Hitchin
fibration has singular fibers.  They remark there that the results in
the thesis can be used in their proof, and further should be useful to
extend the duality over all of $B$.  For a more complete explanation,
we refer to their paper $\cite{DP2}$.

\bigskip

\noindent
{\bf Acknowledgments}\label{ss:thanks} 
Thanks to Tony Pantev for supervising my thesis and his
indispensable help.  Thanks also to Jonathan Block, Ron Donagi, Jim
Stasheff, and Dima Arinkin, and many others, for interesting
conversations and comments on the original draft of the thesis.
\bibliographystyle{amsplain}
\bibliography{bib.bib}

\noindent
email: oren.benbassat@gmail.com

\end{document}